\documentstyle[11pt]{article}
\newcommand{\psn}{{\bf R}^{n} \times {\bf R}^{n}}
\newcommand{\al}{\alpha}

\newcommand{\ovl}{\overline}
\newcommand{\supim}{\mbox{\rm{{\tiny image}}}\,}
\newcommand{\image}{\mbox{\rm{image}}\,}
\newtheorem{definition}{Definition}[section]
\newtheorem{theorem}[definition]{Theorem}
\newtheorem{lemma}[definition]{Lemma}
\newtheorem{corollary}[definition]{Corollary}

\newtheorem{remark}{Remark}
\begin{document}
\baselineskip=20pt
\begin{center}
{\bf {\large  The inverse problem of the Birkhoff-Gustavson
normalization and}}
\\
{\bf {\large {\it ANFER}, Algorithm of Normal Form Expansion and
Restoration}}\footnote{Submitted to J.~Comput. and Appl.~Math..}
\\
\bigskip
Yoshio Uwano\footnote{uwanp@amp.i.kyoto-u.ac.jp}
\\
\medskip
Department of Applied Mathematics and Physics \\
Kyoto University \\
Kyoto 606-8501, Japan \\
\end{center}
\bigskip
\baselineskip=16pt
\begin{abstract}
In the series of papers [1-4], the inverse problem of
the Birkhoff-Gustavson normalization was posed
and studied.
To solve the inverse problem,
the symbolic-computing program named ANFER
(Algorithm of Normal Form Expansion and Restoration)
is written up, with which a new aspect
of the Bertrand and Darboux integrability condition is
found \cite{Uwano2000}.
In this paper, the procedure in ANFER is presented
in mathematical terminology, which is organized on the
basis of the composition of canonical transformations.
\end{abstract}
\baselineskip=12pt
\newpage
\section{Introduction}
It has been recognized that the Birkhoff-Gustavson
normalization \cite{Gustavson} works effectively in
studying various nonlinear dynamical systems.
For example, when a two-degree-of-freedom Hamiltonian
system with a 1:1-resonant equilibrium point is given,
for example, the BG-normalization of its Hamiltonian
around the equilibrium point provides an \lq approximate'
Hamiltonian system whose Hamiltonian is the truncation
of the normalized Hamiltonian up to a finite degree:
The approximate Hamiltonian system provides a good
account of the surface of section with sufficiently
small energies of the given system \cite{Kummer,Cushman}.
Such a good approximation thereby implies that
to find the family of Hamiltonian systems sharing
the same BG-normalization up to a finite-degree
amounts to find a family of Hamiltonian systems
admitting the surface of section similar to each other.
The following question has been posed by the author
as the {\it inverse problem} of the BG-normalization
\cite{Uwano2000,CASC,CPC,JadPhys}:
{\it What kind of polynomial Hamiltonian can be brought into
a given polynomial Hamiltonian in BG-normal form~?}
\par
The Birkhoff-Gustavson normalization to be dealt with this paper is
outlined as follows (cf. \cite{Moser,Uwano2000}):
Let $\psn$ be the phase space
endowed with the Cartesian coordinates $(q,p)$
working as canonical coordinates.
Let $K(q,p)$ be a Hamiltonian function defined on
$\psn$ (or a certain domain of it) which admits
the power-series expansion
\begin{equation}
\label{K}
K(q,p)
=
\sum_{j=1}^{n}\frac{\nu_j}{2}(p_j^2+q_j^2)
+
\sum_{k=3}^{\infty}K_k (q,p),
\end{equation}
around the origin of $\psn$, where each $K_k(q,p)$
($k=3,4,\cdots$) is a homogeneous polynomial of
degree-$k$ in $(q,p)$, and $\{ \nu_j \}$ non-vanishing
constants.
\begin{remark}
\label{convergence}
The convergent radius of the power series (\ref{K})
may vanish \cite{Moser}; this happens to any $K$ that
is not analytic but differentiable around the origin,
for example. In such a case, the power series (\ref{K})
is considered only in a formal sense. We will, however,
often eliminate the word \lq formal' from such formal
power series henceforth.
\end{remark}
The BG-normalization of $K(q,p)$ is made as follows:
Let $(\xi,\eta)$ be another canonical coordinates
of $\psn$ and let the power series
(see Remark~\ref{convergence})
\begin{equation}
\label{W}
W(q,\eta)=\sum_{j=1}^{n}q_j\eta_j 
+
\sum_{k=3}^{\infty} W_{k}(q,\eta)
\end{equation}
be a second-type generating function 
\cite{Goldstein},
a function in the {\it old} position variables $q$ and
the {\it new} momentum variables $\eta$, associated
with the canonical transformation
\begin{equation}
\label{can-W}
(q,p) \rightarrow (\xi,\eta) 
\quad 
\mbox{with}
\quad
p=\frac{\partial S}{\partial q}
\quad \mbox{and} \quad
\xi=\frac{\partial S}{\partial \xi},
\end{equation}
where each $W_{k}(q, \eta)$ ($k=3,4,\cdots $)
is the homogeneous polynomials of degree $k$ 
in $(q, \eta)$. Note that the transformation
(\ref{can-W}) leaves the origin of $\psn$ invariant.
Through the transformation (\ref{can-W}),
the Hamiltonian $K(q,p)$ is brought into
a power series, say $G(\xi, \eta)$, subject to
\begin{equation}
\label{defeq-ord}
G( \frac{\partial W}{\partial \eta}, \eta )
=H ( q , \frac{\partial W}{\partial q}),
\end{equation}
It is easy to see from (\ref{K}), (\ref{W}) and
(\ref{defeq-ord}) that $G(\xi,\eta)$ takes the form
\begin{equation}
\label{G}
G(\xi, \eta)=
\sum_{j=1}^{n}\frac{\nu_j}{2} \left( \eta_j^2 +\xi_j^2
\right) + \sum_{k=3}^{\infty} G_{k}(\xi, \eta),
\end{equation}
where each $G_{k}(\xi,\eta)$ ($k=3,4, \cdots $) is a
homogeneous polynomial of degree-$k$ in $(\xi,\eta)$.
\begin{definition}
\label{def-BGNF}
The power series $G(\xi, \eta)$ is said to be in
the Birkhoff-Gustavson (BG-) normal form
up to degree-$r$ if $G(\xi, \eta)$ satisfies
the Poisson-commuting relation,
\begin{equation}
\label{Poisson}
\left\{ 
\sum_{j=1}^{n}\frac{\nu_j}{2}\left( \eta_j^2 +\xi_j^2
\right) , \;
G_{k}( \xi , \eta) \right\}_{\xi,\eta} =0
\quad (k=3, \cdots , r),
\end{equation}
where $\{ \cdot , \cdot \}_{\xi,\eta}$ denotes
the canonical Poisson bracket to the coordinates
$(\xi,\eta)$ (see \cite{Arnold}).
\end{definition}
The inverse problem of the BG-normalization has been
hence posed as follows \cite{Uwano2000,CASC}:
{\it
For a given power series Hamiltonian $G(\xi,\eta)$
in the BG-normal form (\ref{G}) with 
(\ref{Poisson}, $r=\infty$), identify all the possible
power-series Hamiltonians whish share $G(\xi,\eta)$
as their BG-normalization}. In contrast with the
inverse prblem, we will refer to the problem of
normalizing Hamiltonians into BG-normal form as the
{\it ordinary} problem.
\par
Since elementary algebraic operations, differentiation,
and integration of polynomials have to be repeated many
times to solve the inverse problem, computer algebra
is worth applying to the inverse problem; see
\cite{Uwano2000,CASC,UwanoWEB} for \lq $\mbox{ANFER}$'
and \cite{CPC,JadPhys} for \lq $\mbox{GITA}^{-1}$'.
The procedure in ANFER is based deeply on the composition
of canonical transformation, which is published in
\cite{Uwano2000} only for the two-degree-of-freedom in 1:1
resonance case.
\par
The aim of this paper is to present the procedure in
ANFER in mathematical terminology. After the procedure, 
a new aspect of the Bertrand-Darboux integrability found
through the inverse problem of certain perturbed oscillator
Hamiltonians \cite{Uwano2000} is presented briefly as an
application of ANFER. The contents of this paper is
organized as follows.
\par
Section~2 sets up the inverse problem of the
BG-normalization: For a better understanding of the
inverse problem, the solution of the ordinary problem
is given first. After that, the inverse problem is posed
and solved. Section~3 is devoted to studying the
composition of canonical transformations which provides
a key of organizing ANFER. Those who have taken a course
of analytical mechanics might think that the composition
is well-known already. However, it seems that almost all
the things known of are on the composition of
{\it infinitesimal} transformations. Since the
transformations dealt with here are not infinitesimal,
section~3 is indeed important to organize ANFER.
In section~4, the procedure in ANFER is described.
Section~5 is for an application of ANFER: a new
aspect of the Bertrand-Darboux integrability condition
for certain the perturbed harmonic oscillators \cite{Uwano2000}
is presented. Section~6 is for the concluding remarks.
\section{The inverse problem of the BG-normalization}
The major part of this section is devoted to the inverse
problem of the BG-normalization posed in
\cite{Uwano2000,CASC,CPC,JadPhys}.
Before the inverse problem, however,
we show the way solve the ordinary problem, which will
promote a better understanding of the inverse problem.
\subsection{Solving the ordinary problem}
\quad
Let us start with equating the homogeneous-polynomial part
of degree-$k$ ($k=3,4, \cdots $) in (\ref{defeq-ord}).
Then equation (\ref{defeq-ord}) is put into the series
of equations,
\begin{equation}
\label{eq-ord}
G_{k} (q , \eta ) + (D_{q,\eta} W_{k} ) 
=
K_{k}(q,\eta) + \Phi_k (q,\eta)   \quad (k=3,4, \cdots ),
\end{equation}
where $D_{q,\eta}$ is the differential operator,
\begin{equation}
\label{def-D}
D_{q,\eta}
=
\sum_{j=1}^{n} \nu_j \left( 
q_{j} \frac{\partial}{\partial \eta_{j}}
-
\eta_{j} \frac{\partial}{\partial q_{j}}
\right) .
\end{equation}
The $\Phi_k (q, \eta)$ in (\ref{eq-ord})
is the homogeneous polynomial of degree-$k$
in $(q, \eta)$ which is uniquely determined by
$W_{3}, \cdots , W_{k-1}$, $K_{3}, \cdots , K_{k-1}$,  
$G_{3}, \cdots , G_{k-1}$ given:
In particular, we have $\Phi_3 (q,\eta)=0$ and
\begin{equation}
\label{Phi}
\Phi_4 (q,\eta)=
\sum_{j=1}^{n} \left (
\frac{1}{2}
\left(
\frac{\partial W_3}{\partial q_j} 
\right)^2
+ 
\left.
\frac{\partial K_3}{\partial p_j}
\right\vert_{(q,\eta)}
\frac{\partial W_3}{\partial q_j}
-
\frac{1}{2}
\left(
\frac{\partial W_3}{\partial \eta_j} 
\right)^2
-
\left.
\frac{\partial G_3}{\partial \xi_j}
\right\vert_{(q,\eta)}
\frac{\partial W_3}{\partial \eta_j}
\right) ,
\end{equation}
while $\Phi_k$ with larger $k$ would be of very
complicated form to be described.
\par
To solve equation (\ref{eq-ord}),
the direct-sum decomposition induced by $D_{q,\eta}$
of the spaces of homogeneous polynomials are of great
use. Let us denote by $V_{k}  (q, \eta)$ the vector
space of homogeneous polynomials of degree-$k$ in 
$(q,\eta)$ with real-valued coefficients
($k = 0,1, \cdots $). Since the differential operator
$D_{q,\eta}$ acts linearly on each $V_k (q,\eta)$,
the action of $D_{q,\eta}$ naturally induces the
direct-sum decomposition,
\begin{equation}
\label{decomp}
V_{k} (q,\eta)
=
\image D_{q,\eta}^{(k)}  \oplus \ker D_{q,\eta}^{(k)}
\quad (k=0,1, \cdots),
\end{equation}
of $V_k (q,\eta)$,
where $D_{q,\eta}^{(k)}$ denotes the restriction,
\begin{equation}
\label{restriction}
D_{q,\eta}^{(k)}
=
\left. D_{q,\eta}\right\vert_{V_{k}(q,\eta)}
\quad (k=0,1, \cdots).
\end{equation}
For $D_{q,\eta}$ and $D_{q,\eta}^{(k)}$ ($k=3,4,\cdots$),
we have the following easy to be shown.
\begin{lemma}
\label{lemma-Poisson}
Equation (\ref{Poisson}) is equivalent to
\begin{equation}
\label{Lie-deriv}
(D_{q,\eta}(G_k\vert_{(q,\eta)}))(q, \eta)
=
(D_{q,\eta}^{(k)}(G_k\vert_{(q,\eta)}))(q, \eta)
=0
\quad (k=3,\cdots ,r).
\end{equation}
Namely,
\begin{equation}
\label{Gk-in-ker}
G_k\vert_{(q,\eta)} \in \ker D_{q,\eta}^{(k)}
\quad
(k=3, \cdots , r).
\end{equation}
\end{lemma}
\begin{lemma}
\label{lemma-kernel-D}
Let the coefficients, $\{ \nu_j \}_{j=1, \cdots ,n}$,
of $\sum_{j=1}^n (\nu_j/2)(p_j^2+q_j^2)$ in $K(q,p)$
be said to be independent over ${\bf Z}$ (integers)
if and only if
\begin{equation}
\label{independent}
\sum_{j=1}^n \gamma_j \nu_j =0 \quad
(\gamma_j \in {\bf Z})
\quad \Leftrightarrow \quad
\gamma_j=0 \; (j=1, \cdots , n).
\end{equation}
holds true.
Then, $\ker D_{q,\eta}$ is spanned by the even degree
polynomials,
\begin{equation}
P_{m}=\prod_{j=1}^n (\eta_j^2 +q_j^2)^{m_j}
\qquad
(m_j : \mbox{non-negative integer}, \, j=1,\cdots ,n),
\end{equation}
if and only if $\{ \nu_j \}$ are independent over
${\bf Z}$ (cf. \cite{Moser}). In particular, 
if $\{ \nu_j \}$ are independent over ${\bf Z}$
then $\ker D_{q,\eta}^{(k)} = \{ 0 \}$ for every
odd $k$.
\end{lemma}
We proceed to solving (\ref{eq-ord}) now.
According to (\ref{decomp}), let us decompose
$K_{k}(q,\eta)$ and $\Phi_k(q,\eta)$
($k=3,4,\cdots$) to be
\begin{equation}
\label{K-Phi-decomp}
\begin{array}{l}
\displaystyle{
K_k (q,\eta)
=
K_k^{\supim} (q,\eta) + K_k^{\ker} (q,\eta),
}
\\ \noalign{\vskip 6pt}
\displaystyle{
\Phi_k (q,\eta)
=
\Phi_k^{\supim} (q,\eta) + \Phi_k^{\ker} (q,\eta),
}
\end{array}
\end{equation}
where
\begin{equation}
\label{K-Phi-component}
\begin{array}{l}
\displaystyle{
K_k^{\supim} (q,\eta), \Phi_k^{\supim} (q,\eta) 
\in \image D_{q,\eta}^{(k)},
}
\\ \noalign{\vskip 6pt}
\displaystyle{
K_k^{\ker} (q,\eta), \Phi_k^{\ker} (q,\eta) 
\in \ker D_{q,\eta}^{(k)},
.
}
\end{array}
\end{equation}
Since $G_{k}(q,\eta) \in \ker D_{q,\eta}^{(k)}$
by Lemma \ref{lemma-Poisson} and since
$D_{q,\eta}W_{k} \in \image D_{q,\eta}^{(k)}$ ,
we obtain
\begin{equation}
G_{k}(q,\eta)
=
K_k^{\ker} (q,\eta) +\Phi_k^{\ker} (q,\eta) 
\quad
(k=3,4, \cdots),
\label{solution-ord-G}
\end{equation}
as a solution of (\ref{eq-ord}), 
where $W_k$ is chosen to be
\begin{equation}
\label{solution-ord-W}
W_{k}(q,\eta)
=
\left(
\left. D_{q,\eta}^{(k)}\right\vert
_{\supim D_{q,\eta}^{(k)}}^{-1}
(K_k^{\supim}\vert_{(q,\eta)}
 +\Phi_k^{\supim} )\right)
(q,\eta).
\end{equation}
What is crucial of (\ref{solution-ord-W})
is that $W_{k}(q,\eta) \in \image D_{q,\eta}^{(k)}$
($k=3,4,\cdots$), which ensures the uniqueness
of $G(\xi,\eta)$: For a certain integer
$\kappa \geq 3$, let us consider the sum,
$
{\tilde W}_{\kappa}=
W_{\kappa}+(\mbox{any polynomial in
$\ker D_{q,\eta}^{(\kappa)}$})
$,
where $G_{k}$ and $W_{k}$ with $k < \kappa$ are
given by (\ref{solution-ord-G},\ref{solution-ord-W}).
Even after such a modification, ${\tilde W}_{\kappa}$
satisfy (\ref{solution-ord-W}) still, which will
lead another series of solutions of (\ref{eq-ord})
with $k > \kappa$. Therefore, under the restriction,
$W_k \in \image D_{q,\eta}^{(k)}$ ($k=3,4,\cdots$),
(\ref{solution-ord-G}) with (\ref{solution-ord-W})
is said to be the unique solution of (\ref{defeq-ord}).
To summarize, the ordinary problem is defined
as follows:
\begin{definition}[The ordinary problem]
\label{ord-problem}
For a given Hamiltonian $K(q,p)$ in power series
(\ref{K}), bring $K(q,p)$ into the BG-normal form
$G(\xi,\eta)$ in power series (\ref{G}) which satisfy
(\ref{defeq-ord}) and (\ref{Poisson}), where the
second-type generating function $W$ of the form
(\ref{W}) is chosen to satisfy (\ref{defeq-ord}) and
\begin{equation}
\label{cond-W}
W_{k}(q,\eta) \in \image D_{q,\eta}^{(k)}
\quad (k=3,4,\cdots).
\end{equation}
\end{definition}
\begin{theorem}
\label{theorem-ord}
The BG-normal form $G(\xi,\eta)$ for the Hamiltonian
$K(q,p)$ is given by (\ref{G}) with
(\ref{solution-ord-G}), where the second-type generating
function $W(q,\eta)$ in power series (\ref{W}) is chosen
to satisfy (\ref{solution-ord-W}) and (\ref{cond-W}).
\end{theorem}
\subsection{The inverse problem}
To pose the inverse problem appropriately,
we start with looking the key equation, (\ref{defeq-ord}),
of the ordinary problem into more detail
from a viewpoint of canonical transformations.
With $W(q,\eta)$, let us associate the inverse canonical
transformation,
\begin{equation}
\label{canonical-xi-eta}
(\xi , \eta ) \rightarrow (q,p)
\qquad
\mbox{with}
\qquad
\xi = -\frac{\partial (-W)}{\partial \eta}
\quad
\mbox{and}
\quad  
p = -\frac{\partial (-W)}{\partial q},
\end{equation}
of (\ref{can-W}),
so that $-W(q,\eta)$ is regarded as a
{\it third-type} generating function (Goldstein 1950), 
a function of the {\it old} momentum variables $\eta$
and the {\it new} position variables $q$.
Equation (\ref{defeq-ord}) rewritten as
\begin{equation}
\label{eq-restore}
K( q , -\frac{\partial (-W)}{\partial q} )
=
G ( -\frac{\partial (-W)}{\partial \eta }, \eta) ,
\end{equation}
is then combined with (\ref{canonical-xi-eta})
to show the following.
\begin{lemma}
\label{lemma-restore}
Let $G(\xi,\eta)$ of (\ref{G}) be the BG-normal form
for the Hamiltonian $K(q,p)$ of (\ref{K}), which
satisfies (\ref{defeq-ord}) with a second-type generating
function $W(q, \eta)$. The Hamiltonian $K(q,p)$ is
restored from $G(\xi,\eta)$ through the canonical
transformation (\ref{canonical-xi-eta}) associated with
the third-type generating function $-W(q, \eta)$.
\end{lemma}
Now we can pose the inverse problem in the following way: 
Let the Hamiltonian $H(q,p)$ be written in the form,
\begin{equation}
\label{H}
H(q,p)=\frac{1}{2}\sum_{j=1}^{n}
\left( p_j^2 +q_j^2 \right) 
+ \sum_{k=3}^{\infty} H_{k}(q,p),
\end{equation}
where each $H_{k}(q,p)$ ($k=3,4, \cdots $)
is a homogeneous polynomial of degree-$k$
in $(q,p)$.
Further, let a third-type generating function
$S(q , \eta )$ be written in the form,
\begin{equation}
\label{S}
S(q, \eta )
=- \sum_{j=1}^{n} q_j \eta_j  
- \sum_{k=3}^{\infty} S_{k}(q, \eta) ,
\end{equation}
where each $S_{k}(q,p)$ ($k=3,4, \cdots $) is
a homogeneous polynomial of degree-$k$
in $(q,p)$.
\begin{definition}[The inverse problem]
\label{inv-problem}
For a given BG-normal form, $G(\xi, \eta)$, in power
series (\ref{G}), identify all the Hamiltonians
$H(q,p)$ in power series (\ref{H}) which satisfy
\begin{equation}
\label{defeq-inv}
H( q , -\frac{\partial S}{\partial q} )
=G ( -\frac{\partial S}{\partial \eta }, \eta )
,
\end{equation}
where the third-type generating function $S(q,\eta)$
in power series (\ref{S}) is chosen to satisfy
(\ref{defeq-inv}) and
\begin{equation}
\label{cond-S}
S_k(q,\eta) \in \image D_{q,\eta}^{(k)} 
\quad (k=3,4,\cdots).
\end{equation}
\end{definition}
\subsection{Solving the inverse problem}
We solve the inverse problem in the following way.
On equating the homogeneous-polynomial part of
degree-$k$ in (\ref{defeq-inv}), equation
(\ref{defeq-inv}) is put into the series of equations,
\begin{equation}
\label{eq-inv}
H_{k}(q, \eta)-(D_{q,\eta}S_{k})(q, \eta)
= G_{k}(q, \eta) - \Psi_{k}(q, \eta)
\qquad (k=3,4, \cdots),
\end{equation}
where $D_{q, \eta}$ is given by (\ref{def-D}).
The $\Psi_{k} (q,\eta)$ is the homogeneous polynomial
of degree-$k$ in $(q,\eta)$ determined uniquely by
$H_{3}$, $\cdots$, $H_{k-1}$,
$G_{3}$, $\cdots$, $G_{k-1}$, $S_{3}$, $\cdots$,
$S_{k-1}$ given. In particular, we have
$\Psi_3 (q,\eta)=0$ and
\begin{equation}
\label{Psi}
\Psi_4 (q,\eta)
=
\sum_{j=1}^{n} \left (
\frac{1}{2}
\left(
\frac{\partial S_3}{\partial q_j} 
\right)^2
+ 
\left.
\frac{\partial H_3}{\partial p_j}
\right\vert_{(q,\eta)}
\frac{\partial S_3}{\partial q_j}
-
\frac{1}{2}
\left(
\frac{\partial S_3}{\partial \eta_j} 
\right)^2
-
\left.
\frac{\partial G_3}{\partial \xi_j}
\right\vert_{(q,\eta)}
\frac{\partial S_3}{\partial \eta_j}
\right)
,
\end{equation}
while $\Psi_k$ with larger $k$ would be of very
complicated form to be described. Like in the
ordinary problem, we solve (\ref{eq-inv}) for
$H_{k}$ and $S_{k}$ by using the direct-sum
decomposition (\ref{decomp}) of $V_k (q,\eta)$,
the vector spaces of homogeneous polynomials of
degree-$k$ ($k=3,4,\cdots$). Let us decompose
$H_k$ and $\Psi_k$ to be
\begin{equation}
\label{H-Psi-decomp}
\begin{array}{l}
\displaystyle{
H_k (q,\eta)
=
H_k^{\supim}(q,\eta) + H_k^{\ker} (q,\eta),
}
\\ \noalign{\vskip 6pt}
\displaystyle{
\Psi_k (q,\eta)
=
\Psi_k^{\supim} (q,\eta) + \Psi_k^{\ker} (q,\eta),
}
\end{array}
\end{equation}
where
\begin{equation}
\label{H-Psi-component}
\begin{array}{l}
\displaystyle{
H_k^{\supim} (q,\eta), \Psi_k^{\supim} (q,\eta)
\in \image D_{q,\eta}^{(k)},
}
\\ \noalign{\vskip 6pt}
\displaystyle{
H_k^{\ker} (q,\eta), \Psi_k^{\ker} (q,\eta)
\in \ker D_{q,\eta}^{(k)},
}
\end{array}
\end{equation}
Then on equating $\ker D_{q,\eta}^{(k)}$-part in
(\ref{eq-inv}), $H_{k}^{\ker}$ is determined to be
\begin{equation}
\label{solution-inv-ker}
H_{k}^{\ker}(q, \eta) 
= G_{k}(q, \eta) - \Psi_{k}^{\ker}(q, \eta).
\end{equation}
Equating $\image D_{q,\eta}$-part of (\ref{eq-inv}),
we have 
\begin{equation}
\label{eq-inv-image}
H_{k}^{\supim}(q, \eta) 
-\left( D_{q, \eta} S_{k} \right)(q, \eta)
= - \Psi_{k}^{\supim}(q, \eta) .
\end{equation}
Since we have the pair of unidentified polynomials,
$H_{k}^{\supim}$ and $S_{k}$ in (\ref{eq-inv-image}),
$H_{k}^{\supim}$ is not determined uniquely in
contrast with $H_{k}^{\ker}$; such a non-uniqueness
is of the very nature of the inverse problem.
Accordingly, equation (\ref{eq-inv-image}) is solved
to as
\begin{eqnarray}
\label{solution-inv-image-H}
&& H_{k}^{\supim}(q,\eta) \in \image D_{q,\eta}^{(k)}:
\; \mbox{chosen arbitrarily},
\\
\label{solution-inv-image-S}
&& S_{k}(q,\eta)
=
\left(
\left. D_{q,\eta}^{(k)}\right\vert
_{\supim D_{q,\eta}^{(k)}}^{-1}
(H_k^{\supim}\vert_{(q,\eta)} +\Psi_k^{\supim} )
\right)(q,\eta),
\end{eqnarray}
($k=3,4,\cdots$). Thus we have the following.
\begin{theorem}
\label{theorem-inv}
For a given BG-normal form $G(\xi,\eta)$ in power
series (\ref{G}), the solution $H(q,p)$ of the
inverse problem is given by (\ref{H}) subject to
(\ref{solution-inv-ker}) and
(\ref{solution-inv-image-H}),
where the third-type generating function
$S(q,\eta)$ in (\ref{defeq-inv}) is chosen to be
(\ref{S}) subject to (\ref{solution-inv-image-S}).
\end{theorem}
\begin{remark}
We have another expression of and $S_{k}$ and
$H^{\supim}_{k}$ equivalent to
(\ref{solution-inv-ker}, \ref{solution-inv-image-H},
\ref{solution-inv-image-S});
\begin{equation}
\label{another-solution-H}
H_{k}^{\supim}(q,\eta)
= \left( D_{q,\eta}^{(k)} S_{k} \right) (q,\eta)
+ \Psi_{k}^{\supim} (q,\eta),
\end{equation}
with
\begin{equation}
\label{another-solution-S}
S_{k}(q, \eta) \in \image D_{q,\eta}^{(k)}: 
\mbox{chosen arbitrarily} .
\end{equation}
Since we usually pay much more interests in $H$
than in $S$, in the inverse problem, it would be
better to take the expression
(\ref{solution-inv-ker}, \ref{solution-inv-image-H},
\ref{solution-inv-image-S}).
\end{remark}
\subsection{The degree-$\rho$ ordinary and inverse
problems}
From a practical point of view, we usually deal with
the BG-normal form Hamiltonians not in power series
but in polynomial. Indeed, as mentioned in section~1,
when we utilize the BG-normalization to provide an
approximate system for a given system, we truncate
the normalized Hamiltonian up to a finite degree.
Hence we naturally come to think of a
\lq finite-degree version' of both the ordinary
and the inverse problems \cite{Uwano2000,CASC}.
\begin{definition}
[The degree-$\rho$ ordinary problem]
\label{deg-rho-ord}
For a given Hamiltonian $K(q,p)$ of the form
(\ref{K}) (possibly in polynomial form) and an
integer $\rho \geq 3$, bring $K(q,p)$ into the
polynomial $G(\xi,\eta)$ of degree-$\rho$
in BG-normal form which satisfy (\ref{defeq-ord})
up to degree-$\rho$, where the second-type generating
function $W(q,\eta)$ in (\ref{defeq-ord}) is chosen
to be the polynomial of degree-$\rho$ subject to
(\ref{cond-W}) with $k=3,\cdots,\rho$.
\end{definition}
\begin{remark}
\label{degree}
In the papers \cite{Uwano2000,CASC}, the degree,
$\rho$, in Definition~\ref{deg-rho-ord} is restricted to be
an even integers, because all the $\nu_j$'s are
set to be $1$ there.
(cf. Lemmas~\ref{lemma-Poisson} and \ref{lemma-kernel-D})
However, since we are dealing with the general $\{ \nu_j \}$
in this paper, the degree $\rho$ has no restriction.
\end{remark}
\begin{definition}[The degree-$\rho$ inverse problem]
\label{deg-rho-inv}
For a given BG-normal form, $G(\xi, \eta)$, of degree-$\rho$
with an integer $\rho \geq 3$, identify all the
polynomial-Hamiltonians $H(q,p)$ of degree-$\rho$
which satisfy (\ref{defeq-inv}) up to degree-$\rho$,
where the third-type generating function $S(q, \eta)$
is chosen to be the polynomial of degree-$\rho$
subject to (\ref{cond-S}) with $k=3,\cdots,\rho$.
\end{definition}
\section{Composition of canonical transformations}
In Section~3, we have shown the way to solve
the inverse problem which is theoretically complete.
From a practical point of view, however,
it is rather difficult to realize the solution
in the present form,
(\ref{solution-inv-ker}, \ref{solution-inv-image-H},
\ref{solution-inv-image-S}), since the calculations
required will be highly combinatorial:
Let us take the 1:1 resonant two-degree-of-freedon case
($n=2$, $\nu_1=\nu_2=2$), for example. In order to
calculate $\Phi_{k}$ (resp. $\Psi_{k}$), the polynomials
$G_{3}, \cdots , G_{k-1}$, $S_{3},\cdots, S_{k-1}$,
and $K_{3}, \cdots , K_{k-1}$
(resp. $H_{3} , \cdots, H_{k-1}$)
have to be kept each of which is in $V_{\ell}$ of
combinatorially rising dimension, 
$
\sum_{h=1}^{4} {4 \choose h}  {\ell-1 \choose h-1}  ,
$ 
where the symbol $ {\cdot \choose \cdot} $ indicates
the binomial coefficient.
\par
To get rid of such a difficulty,
we break the needed calculations into a
series of calculation of simpler form, which is
realized in ANFER (see \cite{Uwano2000,CASC,UwanoWEB}).
This section is hence devoted to the mathematical basis
for organizing ANFER.
\par
A key idea in writing-up ANFER is to realize 
(\ref{solution-inv-ker}, \ref{solution-inv-image-H},
\ref{solution-inv-image-S}) up to degree-$\rho$
by applying not the canonical transformation
associated with $S(q,\eta)$ but the composition
of canonical transformations
\begin{equation}
\label{transf-tau-h}
\tau_{h}:
(\xi^{(h-1)}, \eta^{(h-1)}) \rightarrow 
(\xi^{(h)}, \eta^{(h)})
\quad \mbox{with} \quad (\xi^{(2)}, \eta^{(2)})
=(\xi,\eta)
\quad (r=3,\cdots , \rho),
\end{equation}
associated with the third-type generating
functions,
\begin{equation}
\label{S^h}
S^{(h)}(\xi^{(h)} , \eta^{(h-1)})
=-\sum_{j=1}^{n}\xi^{(h)}_j \eta^{(h-1)}_j
  -S_{h}(\xi^{(h)} , \eta^{(h-1)})
\qquad (h=3, \cdots, \rho),
\end{equation}
where $S_h (\xi^{(h)}, \eta^{(h-1)})$ is the
homogeneous polynomial part of degree-$h$ of
the generating function $S(q,\eta)$ (see (\ref{S}))
with $(\xi^{(h)}, \eta^{(h-1)})$ in place of
$(q,\eta)$.
\par
We start by showing the following on the generating
function associated with the composition of a pair
of canonical transformations. It should be pointed
out that little is known explicitly of the generating
function of the composition of non-infinitesimal
canonical transformations like the following lemma,
while well-known is that of infinitesimal ones (see
\cite{Goldstein}, for example).
\par
\begin{lemma}
\label{lemma-comp-pair}
Let the function $f^{(h)}(u^{(h)},v^{(h-1)})$
($h=1,2$) be the third-type generating functions of
canonical transformations,
\begin{equation}
\label{sigma}
\sigma_{h}: (u^{(h-1)},v^{(h-1)})
           \rightarrow (u^{(h)},v^{(h)})
\quad
(h=1,2)
,
\end{equation}
with
\begin{equation}
\label{can-rel-uv}
u^{(h-1)}=-\frac{\partial f^{(h)}}{\partial v^{(h-1)}},
\quad
v^{(h)}= - \frac{\partial f^{(h)}}{\partial u^{(h)}}
\quad 
(h=1,2),
\end{equation}
each of which leaves the origin of $\psn$ invariant.
Then the composition,
$\sigma_{2} \circ \sigma_{1}:
(u^{(0)},v^{(0)})
           \rightarrow (u^{(2)},v^{(2)})
$,
is generated by the third-type generating function,
\begin{equation}
\label{f-1,2}
f_{1,2}(u^{(2)}, v^{(0)})
=
\sum_{j=1}^n 
{\tilde u}_j {\tilde v}_j
+f^{(1)}({\tilde u}, v^{(0)})
+f^{(2)}(u^{(2)}, {\tilde v}),
\end{equation}
where ${\tilde u}^{(1)}$ and ${\tilde v}^{(1)}$
on the rhs of (\ref{f-1,2}) are the
$C^{\infty}$-functions of $(u^{(2)},v^{(0)})$
uniquely determined around $(u^{(2)},v^{(0)})=(0,0)$
to satisfy
\begin{equation}
\label{implicit-rel}
{\tilde u}=
-
\frac{\partial f^{(2)}}{\partial v^{(1)}}
(u^{(2)}, {\tilde v}),
\qquad
{\tilde v}=
-\frac{\partial f^{(1)}}{\partial u^{(1)}}
({\tilde u} , v^{(0)}).
\end{equation}
\end{lemma}
\noindent
{\bf Proof}:\quad
Since the canonical relation (\ref{can-rel-uv})
yields the identities \cite{Goldstein}
\begin{equation}
\label{u-v-f}
-\sum_{j=1}^{n}
(u^{(h-1)}_j dv^{(h-1)}_j
+v^{(h)}_j du^{(h)}_j )
=d(f^{(h)}(u^{(h)}, v^{(h-1)}))
\quad
(h=1,2),
\end{equation}
we have
\begin{equation}
\label{comp-u-v-f}
-\sum_{j=1}^n
(u^{(0)}_j dv^{(0)}_j + v^{(2)}_j du^{(2)}_j )
=
d
\left(
 \sum_{j=1}^{n} u^{(1)}_j v^{(1)}_j 
   + f^{(1)}(u^{(1)}, v^{(0)}) 
     + f^{(2)}(u^{(2)}, v^{(1)})
\right)
.
\end{equation}
Further, the canonical relation (\ref{can-rel-uv})
restricts $\{ (u^{(h)}, v^{(h)}) \}_{h=0,1,2}$
to the inverse image, $\sigma^{-1}(0,0,0,0)$,
of the $C^{\infty}$-map
\begin{equation}
\label{map-implicit}
\begin{array}{l}
\displaystyle{
\sigma:
(u^{(0)}, v^{(0)}, u^{(1)}, v^{(1)}, u^{(2)}, v^{(2)})
\in {\bf R}^6
}
\\ \noalign{\vskip 6pt}
\displaystyle{
\mapsto
\left(
u^{(0)}+\frac{\partial f^{(1)}}{\partial v^{(0)}},
u^{(1)}+\frac{\partial f^{(2)}}{\partial v^{(1)}},
v^{(1)}+\frac{\partial f^{(1)}}{\partial u^{(1)}},
v^{(2)}+\frac{\partial f^{(2)}}{\partial u^{(2)}}
\right) \in {\bf R}^4 .
}
\end{array}
\end{equation}
Applying the implicit function theorem \cite{Spivak}
to the map $\sigma$, we obtain uniquely the functions,
${\tilde u}$, ${\tilde v}$, ${\widehat u}$ and
${\widehat v}$, of $(u^{(2)}, v^{(0)})$
subject to
$
\sigma ({\widehat u}, v^{(0)}, {\tilde u}, {\tilde v},
u^{(2)}, {\widehat v}) =0
$
around $(0,0)$. The differentiability comes from that
of $\sigma$ due to the implicit function theorem.
The substitution of $({\tilde u}, {\tilde v})$
into $(u^{(1)},v^{(1)})$ in (\ref{comp-u-v-f}) thereby
shows (\ref{f-1,2}) with (\ref{implicit-rel}).
This ends the proof.
\par
We are now in a position to compose the series of
canonical transformations
$\{ \tau_{h} \}_{h=3, \cdots , \rho}$,
given by (\ref{transf-tau-h}) and (\ref{S^h})
by using Lemma~\ref{lemma-comp-pair}.
We show the following for the generating function
of the composition,
$\tau_{\rho} \circ \cdots \circ \tau_3$.
\begin{lemma}
\label{lemma-comp-tau}
For any fixed integer $\rho \geq 4$, the composition,
$\tau_{\rho} \circ \cdots \circ \tau_{3}$,
of the canonical transformations
$\{ \tau_{h} \}_{h=3, \cdots, \rho}$ given by
\begin{equation}
\label{def-tau-h}
\xi^{(h-1)}
=
-\frac{\partial S^{(h)}}{\partial \eta^{(h-1)}}
(\xi^{(h)}, \eta^{(h-1)}),
\quad
\eta^{(h)}
=
-\frac{\partial S^{(h)}}{\partial \xi^{(h)}}
(\xi^{(h)}, \eta^{(h-1)})
\quad (h=3, 4, \cdots ).
\end{equation}
with (\ref{S^h}) and (\ref{transf-tau-h})
is associated with the third-type-3 generating
function ${\cal S}^{(\rho)}(\xi^{(\rho)}, \eta^{(2)})$
defined recursively by
\begin{eqnarray}
\label{cal-S-h=3}
&\phantom{=}&
{\cal S}^{(3)}(\xi^{(3)},\eta^{(2)})
=S^{(3)}(\xi^{(3)},\eta^{(2)}),
\\
\nonumber
&\phantom{=}&
{\cal S}^{(h)}(\xi^{(h)}, \eta^{(2)})
\\
\label{cal-S}
&=&
\sum_{j=1}^{n} 
    {\tilde \xi}^{(h-1)}_j {\tilde \eta}^{(h-1)}_j 
+ {\cal S}^{(h-1)}({\tilde \xi}^{(h-1)}, \eta^{(2})
+ S^{(h)}(\xi^{(h)}, {\tilde \eta}^{(h-1)})
\quad (h=4, \cdots , \rho),
\end{eqnarray}
where $S^{(h)}$ are the third-type generating functions
given by (\ref{S^h}), and ${\tilde \xi}^{(h-1)}$ and
${\tilde \eta}^{(h-1)}$ the $C^{\infty}$-functions of
$(\xi^{(h)}, \eta^{(2)})$ determined uniquely
to satisfy
\begin{equation}
\label{tilde-x-y}
{\tilde \xi}^{(h-1)}
=
-\frac{\partial S^{(h)}}{\partial \eta^{(h-1)}}
(\xi^{(h)}, {\tilde \eta}^{(h-1)}),
\quad
{\tilde \eta}^{(h-1)}
=
-\frac{\partial {\cal S}^{(h-1)}}{\partial \xi^{(h-1)}}
({\tilde \xi}^{(h-1)}, \eta^{(2)})
\quad (h=3, 4, \cdots ),
\end{equation}
around $(\xi^{(h)},\eta^{(2)})=(0,0)$.
\end{lemma}
{\bf Proof}: \quad
We prove the lemma by induction:
The starting case, (\ref{cal-S}) with $h=4$, immediately
follows due to Lemma~\ref{lemma-comp-pair}.
We move on to show the case of $h =r (\geq 5)$ in turn
under the assumption that (\ref{cal-S})
and (\ref{tilde-x-y}) hold true for $h=r-1$.
Applying Lemma~\ref{lemma-comp-pair} to the pair
of transformations, $\tau_r$ with $S^{(r)}$ and
$\tau_{r-1} \circ \cdots \circ \tau_3$
with ${\cal S}^{(r-1)}$, we obtain
${\cal S}^{(r)}(\xi^{(r)},\eta^{(2)})$
given by (\ref{cal-S}) and (\ref{tilde-x-y}) with
$h=r$ as the third-type generating function of
$\tau_r \circ (\tau_{r-1} \circ \cdots \circ \tau_{3})$.
The differentiability of ${\tilde \xi}^{(r-1)}$
and ${\tilde \eta}^{(r-1)}$ follows from that of
the generating functions ${\cal S}^{(r-1)}$ and
$S^{(r)}$ (see the proof of Lemma~\ref{lemma-comp-pair}).
This ends the proof.
\par
Expanding ${\cal S}^{(h)}$ in terms of
$(\xi^{(h)},\eta^{(2)})$,
we have the following, a key theorem of ANFER.
\begin{theorem}
\label{theorem-expansion}
The third-type generating functions,
${\cal S}^{(h)}(\xi^{(h)},\eta^{(2)})$ ($h=4, 5, \cdots$),
admit the power series expansion,
\begin{equation}
\label{expansion}
\begin{array}{l}
\displaystyle{
{\cal S}^{(h)}(\xi^{(h)},\eta^{(2)})
=
-\sum_{j=1}^{n}\xi^{(h)}_{j}\eta^{(2)}_{j} 
- 
\sum_{k=3}^{h} S_{k}(\xi^{(h)},\eta^{(2)}) 
+
o_{h}(\xi^{(h)},\eta^{(2)}),
}
\\ \noalign{\vskip 9pt}
\displaystyle{
\frac{o_{h}(\xi^{(h)},\eta^{(2)})}
{(\sum_{j=1}^{n}({\xi^{(h)}_j}^2
   +{\eta^{(2)}_j}^2))^{h/2}}
 \rightarrow 0
\qquad ((\sum_{j=1}^{n}({\xi^{(h)}_j}^2 
   +{\eta^{(2)}_j}^2))^{1/2}
         \rightarrow 0 ),
}
\end{array}
\end{equation}
where $S_{k}(\xi^{(h)},\eta^{(2)})$ are the homogeneous
part of degree-$k$ in $S(q,\eta)$ (see (\ref{S})) with
$(\xi^{(h)}, \eta^{(2)})$ in place of $(q,\eta)$.
\end{theorem}
\noindent
{\bf Proof}:
We show (\ref{expansion}) by induction. Expanding
(\ref{tilde-x-y}) with $h=4$ in terms of
$(\xi^{(4)},\eta^{(2)})$, we have
\begin{eqnarray}
\label{expand-tilde-x-h=4}
{\tilde \xi}^{(3)}
=\xi^{(4)}
  +\frac{\partial S_{4}}{\partial \eta^{(3)}}
        ({\xi}^{(4)}, {\tilde \eta}^{(3)})
,
\\
\label{expand-tilde-y-h=4}
{\tilde \eta}^{(3)}
=\eta^{(2)}
  +\frac{\partial S_{3}}{\partial \xi^{(3)}}
        ({\tilde \xi}^{(3)}, \eta^{(2)}).
\end{eqnarray}
Since the second term on the rhs of
(\ref{expand-tilde-x-h=4}) is of degree higher
than two and since that of (\ref{expand-tilde-y-h=4})
are of degree higher than one,
(\ref{expand-tilde-x-h=4}) and (\ref{expand-tilde-y-h=4})
are put together to provide a further expansion,
\begin{equation}
\label{expand-tilde-x-y-h=4}
\begin{array}{l}
\displaystyle{
{\tilde \xi}^{(3)}
=\xi^{(4)}
  +\frac{\partial S_{4}}{\partial \eta^{(3)}}
        (\xi^{(4)}, \eta^{(2)})
   +o_{3}^{\prime}
       (\xi^{(4)}, \eta^{(2)})
}
,
\\ \noalign{\vskip 9pt}
\displaystyle{
{\tilde \eta}^{(3)}
=\eta^{(2)}
  +\frac{\partial S_{3}}{\partial \xi^{(3)}}
        (\xi^{(4)}, \eta^{(2)})
   +o_{3}^{\prime \prime}
          (\xi^{(4)}, \eta^{(2)}),
}
\end{array}
\end{equation}
where $o_{3}^{\prime}$ and $o_{3}^{\prime \prime}$
indicate the terms in
$
(\sqrt{\sum_{j=1}^n 
({\xi^{(4)}_j}^2 + {\eta^{(2)}_j}^2)})
$
of degree higher than $three$. Equations
(\ref{S^h},\ref{cal-S},\ref{expand-tilde-x-y-h=4})
are then put together to yield
\begin{eqnarray}
&&
{\cal S}^{(4)}(\xi^{(4)},\eta^{(2)})
\nonumber
\\
&=&
\sum_{j=1}^n
\left(
\xi^{(4)}_j
\left.  +\frac{\partial S_{4}}
    {\partial \eta^{(3)}_j}\right\vert
_{(\xi^{(4)}, \eta^{(2)})} +o^{\prime}_3
\right)
\left(
\eta^{(2)}_j
\left.  +\frac{\partial S_{3}}
   {\partial \xi^{(3)}_j}\right\vert
_{(\xi^{(4)}, \eta^{(2)})} +o^{\prime \prime}_3
\right)
\nonumber
\\
&\phantom{=}&  
- \left\{
\sum_{j=1}^n
\left(
\xi^{(4)}_j
\left. +\frac{\partial S_{4}}
   {\partial \eta^{(3)}_j}\right\vert
_{( \xi^{(4)}, \eta^{(2)})} + \o^{\prime}_3
\right)
\eta^{(2)}_j
+
S_3
 (
  \xi^{(4)}
   \left.  +\frac{\partial S_{4}}{\partial \eta^{(3)}}\right\vert
    _{(\xi^{(4)}, \eta^{(2)})} + o^{\prime}_3 ,
  \eta^{(2)}_j
)
\right\}
\nonumber
\\
& \phantom{=} &
- \left\{
\sum_{j=1}^n \xi^{(4)}_j 
\left(
\eta^{(2)}_j
\left.  +\frac{\partial S_{3}}
   {\partial \xi^{(3)}_j}\right\vert
_{(\xi^{(4)}, \eta^{(2)})} +o^{\prime \prime}_3
\right)
+
S_{4}(\xi^{(4)}, 
  \eta^{(2)}
   \left.  +\frac{\partial S_{3}}
     {\partial \xi^{(3)}}
    \right\vert_{(\xi^{(4)}, \eta^{(2)})} 
      +o^{\prime \prime}_3
)
\right\}
\nonumber
\\
&=&
-\sum_{j=1}^{n}\xi^{(4)}_{j}\eta^{(2)}_{j} 
- 
\sum_{k=3}^{4} S_{k}(\xi^{(4)},\eta^{(2)}) 
+
o_{4}(\xi^{(4)},\eta^{(2)}),
\label{expand-calc-h=4}
\end{eqnarray}
which shows (\ref{expansion}) with $h=4$.
We move on to show (\ref{expansion}) with
$h=r \geq 5$ in turn under the assumption
that (\ref{expansion}) with $h=r-1$ holds
true. In a similar way to get
(\ref{expand-tilde-x-y-h=4}), we obtain
\begin{equation}
\label{expand-tilde-x-y}
\begin{array}{l}
\displaystyle{
{\tilde \xi}^{(r-1)}
=\xi^{(r)}
 +\frac{\partial S_{r}}
        {\partial \eta^{(r-1)}}
        (\xi^{(r)}, \eta^{(2)})
   +o_{r-1}^{\prime}
       (\xi^{(r)}, \eta^{(2)})
}
,
\\ \noalign{\vskip 9pt}
\displaystyle{
{\tilde \eta}^{(r-1)}
=\eta^{(2)}
  + \sum_{k=3}^{r-1} 
   \frac{\partial S_k}{\partial \xi^{(k)}}
        (\xi^{(r)}, \eta^{(2)})
   +o_{r-1}^{\prime \prime}
          (\xi^{(r)}, \eta^{(2)}),
}
\end{array}
\end{equation}
where $o_{r-1}^{\prime}$ and $o_{r-1}^{\prime \prime}$
indicate the terms in
$
(\sqrt{\sum_{j=1}^n 
({\xi^{(r)}_j}^2 + {\eta^{(2)}_j}^2)})
$
of degree higher than $r-1$.
Like (\ref{expand-calc-h=4}), equations (\ref{S^h})
with $h=r$, (\ref{cal-S}) with $h=r-1$ and
(\ref{expand-tilde-x-y-h=4}) are thereby put
together to yield
\begin{eqnarray}
&&
{\cal S}^{(r)}(\xi^{(r)},\eta^{(2)})
\nonumber
\\
&=&
\sum_{j=1}^n
\left(
\xi^{(r)}_j
\left.  +\frac{\partial S_{r}}
        {\partial \eta^{(r-1)}_j}\right\vert
_{(\xi^{(r)}, \eta^{(2)})} +o^{\prime}_{r-1}
\right)
\left(
\eta^{(2)}_j +
\sum_{k=3}^{r-1} 
   \left.  +\frac{\partial S_k}
           {\partial \xi^{(k)}_j}\right\vert
_{(\xi^{(r)}, \eta^{(2)})} 
  +o^{\prime \prime}_{r-1}
\right)
\nonumber
\\
& \phantom{=} &
- \left\{
\sum_{j=1}^n
\left(
\xi^{(r)}_j
\left. +\frac{\partial S_{r}}
       {\partial \eta^{(r-1)}_j}\right\vert
_{( \xi^{(r)}, \eta^{(2)})} + o^{\prime}_{r-1}
\right)
\eta^{(2)}_j
\right.
\nonumber
\\
& \phantom{=}&
\phantom{========}
\left.
+
\sum_{k=3}^{r-1}S_k
 (
  \xi^{(r)}
   \left.  +\frac{\partial S_r}
           {\partial \eta^{(3)}}\right\vert
    _{(\xi^{(r)}, \eta^{(2)})} 
     + o^{\prime}_{r-1} ,
  \eta^{(2)}_j
)
\right\}
\nonumber
\\
&\phantom{=} &
- \left\{
\sum_{j=1}^n \xi^{(r)}_j 
\left(
\eta^{(2)}_j +\sum_{k=3}^{r-1}
\left.  \frac{\partial S_k}
       {\partial \xi^{(k)}_j}\right\vert
_{(\xi^{(r)}, \eta^{(2)})} 
    +o^{\prime \prime}_{r-1}
\right)
\right.
\nonumber
\\
& \phantom{=} &
\phantom{========}
\left. 
+
S_{r}(\xi^{(r)}, 
  \eta^{(2)} +\sum_{k=3}^{r-1}
   \left. \frac{\partial S_k}
         {\partial \xi^{(3)}}
    \right\vert_{(\xi^{(r)}, \eta^{(2)})} 
    +o^{\prime \prime}_{r-1}
)
\right\}
\nonumber
\\
& \phantom{=}&
+ {\tilde o}_r (\xi^{(r)}, \eta^{(2)})
\nonumber
\\
&=&
-\sum_{j=1}^{n}\xi^{(r)}_{j}\eta^{(2)}_{j} 
- 
\sum_{k=3}^{r} S_{k}(\xi^{(r)},\eta^{(2)}) 
+
o_{r}(\xi^{(r)},\eta^{(2)}),
\label{expand-calc}
\end{eqnarray}
which proves our assertion.
\begin{corollary}
\label{cal-S=S}
${\cal S}^{(h)}(q,\eta)$ coincides
with $S(q,\eta)$ up to degree-$h$ if expanded.
\end{corollary}
\section
{The procedure in ANFER}
In the previous section, we have studied the
composition of canonical transformations necessary
to organize ANFER. This section is devoted to the
solving procedure in ANFER in which the composition
of transformations is utilized effectively.
\subsection{The solution of the degree-$\rho$
inverse problem}
We show that the solution of the degree-$\rho$
inverse problem is obtained by applying the
composition $\tau_{\rho} \circ \cdots \circ \tau_3 $
of canonical transformations studied in section.~3.
Let us fix a BG-normal form power series
(possibly a polynomial) Hamiltonian $G(\xi,\eta)$
of the form (\ref{G}). For $G(\xi,\eta)$ fixed,
we define the power-series Hamiltonians
$H^{(h)}(\xi^{(h)},\eta^{(h)})$
($h=3, \cdots , \rho$) by the reccurent formula
\begin{equation}
\label{def-H^h}
H^{(h)} \circ \tau_h = H^{(h-1)}
\quad (h=3, 4, \cdots ) \quad
\mbox{with}
\quad
H^{(2)}(\xi, \eta)=G(\xi,\eta).
\end{equation}
Namely, $H^{(h)}$ satisfy the equation
\begin{equation}
\label{H^2-H^h}
(H^{(h)} \circ \tau_h \circ \cdots \circ \tau_3)
(\xi , \eta)
=
H^{(2)} (\xi, \eta)
(=G(\xi, \eta))
\quad (h=3,4, \cdots),
\end{equation}
which is put together with Lemma~\ref{lemma-comp-tau}
to show 
\begin{equation}
\label{H^h-H^2-S}
H^{(h)}(\xi^{(h)},
   -\frac{\partial {\cal S}^{(h)}}{\xi^{(h)}})
=
H^{(2)}
(-\frac{\partial {\cal S}^{(h)}}{\partial \eta^{(2)}},
   \eta^{(2)})
\quad (h=3,4,  \cdots),
\end{equation}
where ${\cal S}^{(h)}(\xi^{(h)}, \eta^{(2)})$ is given by
(\ref{cal-S}) with (\ref{tilde-x-y}).
\par
Let us recall the equation (\ref{eq-inv}) equivalent
to the defining equation (\ref{defeq-inv}) of the
inverse problem. From (\ref{eq-inv}), we see that
only the homogeneous parts in $S$ and $H$ of degree
lower than $h$, and those in $G$ of degree lower than
or equal to $h$ (cf. (\ref{Psi})) are necessary to
determine the homogeneous part $H_h$ of degree-$h$
in $H$ This observation is thereby put together with
Corollary~\ref{cal-S=S} to show the following.
\begin{theorem}
\label{theorem-H^h=H}
For a given BG-normal form power-series $G(\xi,\eta)$,
the power-series Hamiltonian $H^{(h)}(q,p)$ coincides
up to degree-$h$ ($h=3,4,\cdots$) with the solution
of the inverse problem, $H(q,p)$.
\end{theorem}
\begin{corollary}
\label{col^H^h=H}
The solution of the degree-$\rho$ inverse problem is
given by $H^{(\rho)}(q,p)$ truncated up to degree-$\rho$.
\end{corollary}
\subsection{The solving procedure in ANFER}
We are now in a position to present the procedure of
solving equation (\ref{eq-inv}) with $h=3, \cdots \rho$
(equivalently, (\ref{eq-inv}) up to degree-$\rho$)
for the degree-$\rho$ inverse problem.
In ANFER, the series of equations (\ref{def-H^h}) are
solved as follows: To be more precise, what we solve
is
\begin{equation}
\label{H^h-S-recursion}
H^{(h)}
(\xi^{(h)}, 
   -\frac{\partial S^{(h)}}{\partial \xi^{(h)}})
=
H^{(h-1)}
(-\frac{\partial S^{(h)}}{\partial \eta^{(h-1)}},
\eta^{(h-1)})
\quad
(h=3,\cdots , \rho )
\end{equation}
with
\begin{equation}
\label{initial-H^2}
H^{(2)}(\xi^{(2)},\eta^{(2)})
=
G(\xi^{(2)},\eta^{(2)})
\end{equation}
equivalent to (\ref{def-H^h}).
For convenience, we will refer to the stage of
solving (\ref{H^h-S-recursion}) with $h=r$ as
the \lq Stage-$r$' hence force. Before starting
the Stage-three, (\ref{initial-H^2}) is assumed
to have been proceeded already.
Further, it is convenient to express
$H^{(h)}(\xi^{(h)}, \eta^{(h)})$ to be
\begin{equation}
\label{expand-H^h}
H^{(h)}(\xi^{(h)}, \eta^{(h)})
=
\frac{1}{2}\sum_{j=1}^{n}
  ( {\eta_j^{(h)}}^2 + {\xi_j^{(h)}}^2)
+
\sum_{k=3}^{\infty}
H^{(h)}_k ( \xi^{(h)}, \eta^{(h)} )
\end{equation}
like $H(q,p)$ (see (\ref{H})), where each
$H^{(h)}_k ( \xi^{(h)}, \eta^{(h)} )$ ($k=3,4, \cdots$)
is the homogeneous polynomial part of degree-$k$ in
$( \xi^{(h)}, \eta^{(h)} )$.
\par\smallskip
\noindent
[{\bf Stage}-$r$]:
\newline
At the Stage-$r$, (\ref{H^h-S-recursion}) with $h=r$
is solved. Equating the homogeneous part of degree-$k$
($k=3, \cdots , \rho$) in (\ref{H^h-S-recursion})
with $h=r$, we have the series of equation,
\begin{eqnarray}
\label{Hrk-k<r}
&&
H^{(r)}_{k}(\xi^{(r)},\eta^{(r-1)})
=
H^{(r-1)}_{k}(\xi^{(r)},\eta^{(r-1)})
\qquad (k=3, \cdots, r-1),
\\
\label{Hrk-k=r}
&&
H^{(r)}_{r}(\xi^{(r)} , \eta^{(r-1)})
=
\left(
D_{\xi^{(r)}, \eta^{(r-1)}}S_{r}
\right)(\xi^{(r)},\eta^{(r-1)})
+
H^{(r-1)}_{r}(\xi^{(r)},\eta^{(r-1)})
\quad
(k=r),
\\
\label{Hrk-k>r}
&&
H^{(r)}_{k}(\xi^{(r)}, \eta^{(r-1)})
=
H^{(r-1)}_{k}(\xi^{(r)},\eta^{(r-1)})
+
\Theta^{(r)}_{k}(\xi^{(r)}, \eta^{(r-1)})
\qquad (k=r+1, \cdots , \rho).
\end{eqnarray}
Here $\Theta^{(r)}_{k}(\xi^{(r)}, \eta^{(r-1)})$
$(k=r+1, \cdots , \rho)$ in (\ref{Hrk-k>r})
are the  homogeneous polynomials of degree-$k$
given by
\begin{eqnarray}
\nonumber
&&
\Theta^{(r)}_{k}(\xi^{(r)},\eta^{(r-1)})
\\
\label{Thetark}
&=& \sum_{ \vert {\bf \al} \vert=1 }
^{ \left[ \frac{k-2}{r-2} \right] } 
\frac{1}{{\bf \al}!}
\left[
\left(
\left.
\frac{\partial S_{r}}{\partial \eta}
\right\vert_{(\xi^{(r)}, \eta^{(r-1)})}
\right)^{\al}
\left.
\left(
\left(
\frac{\partial}{\partial \xi^{(r-1)}}
\right)^{\al} 
H^{(r-1)}_{k - (r-2)\vert \al \vert}
\right)
\right\vert_{(\xi^{(r)}, \eta^{(r-1)})}
\right.
\\
\nonumber
&&
\phantom{======}
- \left.
\left(
\left.
\frac{\partial S_{r}}{\partial \xi}
\right\vert_{( \xi^{(r)}, \eta^{(r-1)})}
\right)^{\al}
\left.
\left(
\left(
\frac
{\partial}{\partial \eta^{(r)}}
\right)^{\al} 
H^{(r)}_{k - (r-2)\vert \al \vert}
\right)
\right\vert_{(\xi^{(r)}, \eta^{(r-1)})}
\right] ,
\end{eqnarray}
where $[ (k-2)/(r-2) ]$ stands for the integer-part
of $(k-2)/(r-2)$, and $\al=(\al_1, \cdots , \al_n)$
does the multi-index with nonnegative-integer valued
components associated with the notations,
\begin{equation}
\label{alpha}
\begin{array}{ll}
\displaystyle{
\vert \al \vert = \sum_{j=1}^{n} \al_{j} ,}
& \quad
\displaystyle{
\al ! = \al_{1} ! \, \cdots \al_{n}! ,
} 
\\ \noalign{\vskip 6pt}
\displaystyle{
\left(
\frac{\partial S_{r}}{\partial \xi}
\right)^{\al}
=
\left(
\frac{\partial S_{r}}{\partial \xi_{1}}
\right)^{\al_{1}}
\cdots
\left(
\frac{\partial S_{r}}{\partial \xi_{n}}
\right)^{\al_{n}} ,
}
& \quad
\displaystyle{
\left(
\frac{\partial S_{r}}{\partial \eta}
\right)^{\al}
=
\left(
\frac{\partial S_{r}}{\partial \eta_{1}}
\right)^{\al_{1}}
\cdots
\left(
\frac{\partial S_{r}}{\partial \eta_{n}}
\right)^{\al_{n}},
}
\\ \noalign{\vskip 6pt}
\displaystyle{
\left(
\frac{\partial}{\partial \xi}
\right)^{\al}
=
\left(
\frac{\partial}{\partial \xi_{1}}
\right)^{\al_{1}}
\cdots
\left(
\frac{\partial}{\partial \xi_{n}}
\right)^{\al_{n}} ,
}
& \quad

\displaystyle{
\left(
\frac{\partial}{\partial \eta}
\right)^{\al}
=
\left(
\frac{\partial}{\partial \eta_{1}}
\right)^{\al_{1}}
\cdots
\left(
\frac{\partial}{\partial \eta_{n}}
\right)^{\al_{n}}.
}
\end{array}
\end{equation}
\begin{remark}
We do not have to deal with the homogeneous part in
(\ref{H^h-S-recursion}) of the degree higher than
$\rho$ taking the observation to (\ref{eq-inv})
into account made above Theorem~\ref{theorem-H^h=H}.
\end{remark}
Since Theorem~\ref{theorem-H^h=H} and (\ref{Hrk-k<r})
are put together to imply that each $H^{(r)}_k$ with
$k < r$ has been identified already as $H^{(k)}_k$ at
the stage-$k$, we do not have to deal with 
(\ref{Hrk-k<r}) a lot: Only the replacement of
$H^{(r-1)}_k$ to $H^{(r)}_k$ ($k < r$) is made.
\par
We solve (\ref{Hrk-k=r}) and (\ref{Hrk-k>r})
in turn by using the decomposition,
\begin{equation}
\label{Hrk-Thetark-decomp}
\begin{array}{l}
\displaystyle{
H^{(r)}_k (\xi^{(r)}, \eta^{(r-1)})
=
{H^{(r)}_k}^{\supim} (\xi^{(r)}, \eta^{(r-1)})
+
{H^{(r)}_k}^{\ker} (\xi^{(r)}, \eta^{(r-1)})
\qquad
(k=r, \cdots, \rho)
}
\\ \noalign{\vskip 6pt}
\displaystyle{
\Theta^{(r)}_k (\xi^{(r)}, \eta^{(r-1)})
=
{\Theta^{(r)}_k}^{\supim} (\xi^{(r)}, \eta^{(r-1)})
+
{\Theta^{(r)}_k}^{\ker} (\xi^{(r)}, \eta^{(r-1)})
\qquad
(k=r+1, \cdots, \rho)
}
\end{array}
\end{equation}
where
\begin{equation}
\label{Hrk-Thetark-component}
\begin{array}{l}
\displaystyle{
{H_k^{(r)}}^{\supim} (\xi^{(r)},\eta^{(r-1)}),
{\Theta_k^{(r)}}^{\supim} (\xi^{(r)}, \eta^{(r-1)}) 
\in \image D_{\xi^{(r)},\eta^{(r-1)}}^{(k)},
}
\\ \noalign{\vskip 6pt}
\displaystyle{
{H_k^{(r)}}^{\ker} (\xi^{(r)}, \eta^{(r-1)}),
{\Theta_k^{(r)}}^{\ker} (\xi^{(r)},\eta^{(r-1)})
\in \ker D_{\xi^{(r)}, \eta^{(r-1)}}^{(k)}.
}
\end{array}
\end{equation}
Like in the way to solve (\ref{eq-inv}) in
subsection~2.3, (\ref{Hrk-k=r}) is solved to be
\begin{eqnarray}
\label{Hrk-k=r-ker}
&&
{H_r^{(r)}}^{\ker}(\xi^{(r)},\eta^{(r-1)})
=
{H_r^{(r-1)}}^{\ker}(\xi^{(r)},\eta^{(r-1)})
\\
\label{Hrk-k=r-image}
&&
{H_r^{(r)}}^{\supim}(\xi^{(r)},\eta^{(r-1)}))
\in \image D_{\xi^{(r)}, \eta^{(r-1)}}^{(r)}
:
\mbox{chosen arbitrarily},
\\
&&
\nonumber
S_{r}(\xi^{(r)}, \eta^{(r-1)})
\\
\label{S_r}
&&
=
\left(
\left. D_{\xi^{(r)},\eta^{(r-1)}}^{(r)}
\right\vert
_{\supim D_{\xi^{(r)},\eta^{(r-1)}}^{(r)}}^{-1}
({H_r^{(r)}}^{\supim}
   - {H_r^{(r-1)}}^{\supim}\vert
_{(\xi^{(r)}, \eta^{(r-1)})})
\right)(\xi^{(r)}, \eta^{(r-1)}).
\end{eqnarray}
We move on to (\ref{Hrk-k>r}) with $k=s \, (>r)$:
Since we have identified $H^{(r)}_k$s with $k < s$
and $S_r$ already before solving (\ref{Hrk-k>r})
with $k=s (>r)$, $\Theta^{(r)}_s$ can be identified
completely by (\ref{Thetark}) with $k=s$.
After the calculation of $\Theta^{(r)}_s$, the
$H^{(r)}_k$ with $k=s (>r)$ is thereby identified
by (\ref{Hrk-k>r}).
\par
In ANFER, the solving procedure described above is
realized in the symbolic computing language, REDUCE3.6
or later. Till to now, a prototype of ANFER for the
1:1 resonant two-degree-of-freedom systems
($n=2$, $\nu_1=\nu_2=1$) has been written up by the
author. The source-code of this prototype is available
on the web page, http://yang.amp.i.kyoto-u.ac.jp/\~{}uwano/,
\cite{UwanoWEB}.
\par
In closing this section, we wish to compare the procedures,
(\ref{eq-inv})-(\ref{solution-inv-image-S})
with (\ref{Hrk-k<r})-(\ref{S_r}).
Through (\ref{eq-inv})-(\ref{solution-inv-image-S})
with $k=r$, $S_3 , \cdots , S_{r-1}$, $G_3 , \cdots , G_{r}$,
and $H_3 , \cdots , H_{r}$ have to be kept.
In contrast with this, $S_{r-1}$,
$H^{(r-1)}_3, \cdots,H^{(r-1)}_{\rho}$, and
$H^{(r)}_3,  \cdots , H^{(r)}_{\rho}$ have to be kept.
Hence, the procedure given in this section would
contributes to the memory-saving.
\section{Application}
Although only a prototype of ANFER exists for a 1:1
resonant case, ANFER has worked very effectively to
find a new deep relation between the Bertrand-Darboux
integrability condition (BDIC) (see
\cite{Marshall,Hietarinta,Whittaker}, for example)
for the perturbed harmonic oscillators with
homogeneous-{\it cubic} potentials (PHOCPs) and for
the oscillators with homogeneous-{\it quartic}
potentials (PHOQPs) \cite{Uwano2000}. The BDIC provides
not only a sufficient condition for the integrability
but also a separability, which will be reviewed briefly
in Appendix.
\subsection{
The degree-$4$ ordinary and inverse problems for PHOCP}
We start by solving the degree-$4$ ordinary problem
for the PHOCP-Hamiltonian
\begin{equation}
\label{cal-F3}
{\cal F}^{(3)}(q,p)
=
\frac{1}{2}\sum_{j=1}^2 (p_j^2 + q_j^2)
+
(f_1q_1^3 + f_2q_1^2 +f_3q_1q_2^2 +f_4 q_2^3),
\end{equation}
where $(f_h)$ are real-valued parameters.
Note that if we choose $(f_h)$ to be
\begin{equation}
\label{f-HH}
f_1=f_3=0, \quad f_2=1, \quad f_4=\mu,
\end{equation}
${\cal F}^{(3)}(q,p)$ becomes the well-known one-parameter
H{\' e}non-Heiles Hamiltonian \cite{Kummer,Cushman}.
\par
Through ANFER (\cite{UwanoWEB}), we obtain the BG-normal form
${\cal G}$ of the following form as the solution of the
degree-$4$ ordinary problem \cite{Uwano2000};
\begin{eqnarray}
&&
{\cal G}(\xi ,\eta)
=\frac{1}{2}\left( \zeta_1{\ovl \zeta}_1 + \zeta_2{\ovl \zeta}_2
\right)
\nonumber
\\
&&
\phantom{{\cal G}(\xi ,\eta)x}
-\frac{15}{16}
(f_1^2\zeta_1^2{\ovl \zeta}_1^2+f_4^2\zeta_2^2{\ovl \zeta}_2^2)
-\frac{3}{4}(f_1f_3\zeta_1\zeta_2{\ovl \zeta}_1{\ovl \zeta}_2
            +f_2f_4\zeta_1\zeta_2{\ovl \zeta}_1{\ovl \zeta}_2)
\nonumber
\\
&&
\phantom{{\cal G}(\xi ,\eta)x}
-\frac{5}{8}(f_1f_2\zeta_1^2{\ovl \zeta}_1{\ovl \zeta}_2
+f_1f_2\zeta_1\zeta_2{\ovl \zeta}_1^2
+f_3f_4\zeta_1\zeta_2{\ovl \zeta}_2^2
+f_3f_4\zeta_2^2{\ovl \zeta}_1{\ovl \zeta}_2)
\nonumber
\\
&&
\phantom{{\cal G}(\xi ,\eta)x}
-\frac{5}{24}(f_2f_3\zeta_1^2{\ovl \zeta}_1{\ovl \zeta}_2
+f_2f_3\zeta_1\zeta_2{\ovl \zeta}_1^2
+f_2f_3\zeta_1\zeta_2{\ovl \zeta}_2^2
+f_2f_3\zeta_2^2{\ovl \zeta}_1{\ovl \zeta}_2)
\label{cal-G}
\\
&&
\phantom{{\cal G}(\xi ,\eta)x}
-\frac{1}{6}(f_2^2\zeta_1\zeta_2{\ovl \zeta}_1{\ovl \zeta}_2
+f_3^2\zeta_1\zeta_2{\ovl \zeta}_1{\ovl \zeta}_2)
-\frac{5}{48}(f_2^2\zeta_1^2{\ovl \zeta}_1^2
+f_3^2\zeta_2^2{\ovl \zeta}_2^2)
\nonumber
\\
&&
\phantom{{\cal G}(\xi ,\eta)x}
-\frac{1}{8}(f_2^2\zeta_1^2{\ovl \zeta}_2^2
+f_2^2\zeta_2^2{\ovl \zeta}_1^2
+f_3^2\zeta_1^2{\ovl \zeta}_2^2
+f_3^2\zeta_2^2{\ovl \zeta}_1^2)
\nonumber
\\
&&
\phantom{{\cal G}(\xi ,\eta)x}
+\frac{1}{16}(f_1f_3\zeta_1^2{\ovl \zeta}_2^2+f_1f_3\zeta_2^2{\ovl \zeta}_1^2
+f_2f_4\zeta_1^2{\ovl \zeta}_2^2+f_2f_4\zeta_2^2{\ovl \zeta}_1^2)
.
\nonumber
\end{eqnarray}
\indent
\par
We solve the degree-$4$ inverse problem
for the BG-normal form ${\cal G}$ in turn:
By ANFER, we have the following polynomial of degree-$4$
as the solution \cite{Uwano2000};
\begin{equation}
\label{cal-H}
{\cal H}(q,p)
=
\frac{1}{2}\sum_{j=1}^{2}
\left( p_j^2 + q_j^2 \right)
+
{\cal H}_3(q,p)
+
{\cal H}_4(q,p)
,
\end{equation}
with
\begin{eqnarray}
&&
\nonumber
{\cal H}_3(q,p)=
a_1z_{1}^3
+a_2z_{1}^2z_{2}+a_3z_{1}z_{2}^2
+a_4z_{2}^3+a_5z_{1}^2{\ovl z}_{1}
+a_6z_{1}^2{\ovl z}_{2}
\\
&&
\nonumber
\phantom{{\cal H}_3(q,p)=}
+a_7z_{1}z_{2}{\ovl z}_{1}+a_8z_{1}z_{2}{\ovl z}_{2}
+a_9z_{2}^2{\ovl z}_{1}+a_{10}z_{2}^2{\ovl z}_{2}
\\
&&
\nonumber
\phantom{{\cal H}_3(q,p)=}
+{\ovl a}_1{\ovl z}_{1}^3
+{\ovl a}_2{\ovl z}_{1}^2{\ovl z}_{2}
+{\ovl a}_3{\ovl z}_{1}{\ovl z}_{2}^2
+{\ovl a}_4{\ovl z}_{2}^3+{\ovl a}_5z_{1}{\ovl z}_{1}^2
+{\ovl a}_6z_{2}{\ovl z}_{1}^2
\\
&&
\label{cal-H3}
\phantom{{\cal H}_3(q,p)=}
+{\ovl a}_7z_{1}{\ovl z}_{1}{\ovl z}_{2}
+{\ovl a}_8z_{2}{\ovl z}_{1}
{\ovl z}_{2}+{\ovl a}_9z_{1}{\ovl z}_{2}^2
+{\ovl a}_{10}z_{2}{\ovl z}_{2}^2 ,
\end{eqnarray}
and
\begin{eqnarray}
{\cal H}_4(q,p)&=&
c_1z_{1}^4+c_2z_{1}^3z_{2}+c_3z_{1}^2z_{2}^2+c_4z_{1}z_{2}^3
\nonumber
\\
&&
+c_5z_{2}^4+c_6z_{1}^3{\ovl z}_{1}
+c_7z_{1}^3{\ovl z}_{2}
+c_8z_{1}^2z_{2}{\ovl z}_{1}
\nonumber
\\
&&
+c_9z_{1}^2z_{2}{\ovl z}_{2}
+c_{10}z_{1}z_{2}^2{\ovl z}_{1}
+c_{11}z_{1}z_{2}^2{\ovl z}_{2}
+c_{12}z_{2}^3{\ovl z}_{1}
+c_{13}z_{2}^3{\ovl z}_{2}
\nonumber
\\
&&
+{\ovl c}_1{\ovl z}_{1}^4
+{\ovl c}_2{\ovl z}_{1}^3{\ovl z}_{2}
+{\ovl c}_3{\ovl z}_{1}^2{\ovl z}_{2}^2
+{\ovl c}_4{\ovl z}_{1}{\ovl z}_{2}^3
\nonumber
\\
&&
+{\ovl c}_5{\ovl z}_{2}^4
+{\ovl c}_6z_{1}{\ovl z}_{1}^3
+{\ovl c}_7z_{2}{\ovl z}_{1}^3
+{\ovl c}_8z_{1}{\ovl z}_{1}^2{\ovl z}_{2}
\nonumber
\\
&&
+{\ovl c}_9z_{2}{\ovl z}_{1}^2{\ovl z}_{2}
+{\ovl c}_{10}z_{1}{\ovl z}_{1}{\ovl z}_{2}^2
+{\ovl c}_{11}z_{2}{\ovl z}_{1}{\ovl z}_{2}^2
+{\ovl c}_{12}z_{1}{\ovl z}_{2}^3
+{\ovl c}_{13}z_{2}{\ovl z}_{2}^3
\nonumber
\\
&&
+8(
a_6{\ovl a}_6z_{1}z_{2}{\ovl z}_{1}{\ovl z}_{2}
+a_9{\ovl a}_9z_{1}z_{2}{\ovl z}_{1}{\ovl z}_{2}
+a_5{\ovl a}_6z_{1}z_{2}{\ovl z}_{1}^2
\nonumber
\\
&&
+a_6{\ovl a}_5z_{1}^2{\ovl z}_{1}{\ovl z}_{2}
+a_9{\ovl a}_{10}z_{2}^2{\ovl z}_{1}{\ovl z}_{2}
+a_{10}{\ovl a}_9z_{1}z_{2}{\ovl z}_{2}^2
)
\nonumber
\\
&&
+6(
a_1{\ovl a}_1z_{1}^2{\ovl z}_{1}^2
+a_4{\ovl a}_4z_{2}^2{\ovl z}_{2}^2
+a_5{\ovl a}_5z_{1}^2{\ovl z}_{1}^2
+a_{10}{\ovl a}_{10}z_{2}^2{\ovl z}_{2}^2
)
\nonumber
\\
&&
+4(
a_1{\ovl a}_2z_{1}^2{\ovl z}_{1}{\ovl z}_{2}
+a_8{\ovl a}_9z_{1}^2{\ovl z}_{2}^2
+a_3{\ovl a}_4z_{1}z_{2}{\ovl z}_{2}^2
+a_5{\ovl a}_8z_{1}z_{2}{\ovl z}_{1}{\ovl z}_{2}
\nonumber
\\
&&
+a_6{\ovl a}_7z_{1}^2{\ovl z}_{2}^2
+a_6{\ovl a}_8z_{1}z_{2}{\ovl z}_{2}^2
+a_7{\ovl a}_9z_{1}^2{\ovl z}_{1}{\ovl z}_{2}
+a_7{\ovl a}_{10}z_{1}z_{2}{\ovl z}_{1}{\ovl z}_{2}
)
\nonumber
\\
\phantom{=}
&&
+4(
a_2{\ovl a}_1z_{1}z_{2}{\ovl z}_{1}^2
+a_4{\ovl a}_3z_{2}^2{\ovl z}_{1}{\ovl z}_{2}
+a_8{\ovl a}_5z_{1}z_{2}{\ovl z}_{1}{\ovl z}_{2}
+a_7{\ovl a}_6z_{2}^2{\ovl z}_{1}^2
\nonumber
\\
&&
+a_8{\ovl a}_6z_{2}^2{\ovl z}_{1}{\ovl z}_{2}
+a_9{\ovl a}_7z_{1}z_{2}{\ovl z}_{1}^2
+a_{10}{\ovl a}_7z_{1}z_{2}{\ovl z}_{1}{\ovl z}_{2}
+a_9{\ovl a}_8z_{2}^2{\ovl z}_{1}^2
)
\nonumber
\\
&&
+\frac{8}{3}(
a_2{\ovl a}_2z_{1}z_{2}{\ovl z}_{1}{\ovl z}_{2}
+a_3{\ovl a}_3z_{1}z_{2}{\ovl z}_{1}{\ovl z}_{2}
)
\nonumber
\\
&&
+2(
-a_6{\ovl a}_6z_{1}^2{\ovl z}_{1}^2
+a_7{\ovl a}_7z_{1}^2{\ovl z}_{1}^2
+a_8{\ovl a}_8z_{2}^2{\ovl z}_{2}^2
-a_9{\ovl a}_9z_{2}^2{\ovl z}_{2}^2
)
\nonumber
\\
&&
+2(
a_1{\ovl a}_3z_{1}^2{\ovl z}_{2}^2
+a_2{\ovl a}_4z_{1}^2{\ovl z}_{2}^2
+a_5{\ovl a}_7z_{1}^2{\ovl z}_{1}{\ovl z}_{2}
-a_5{\ovl a}_9z_{1}^2{\ovl z}_{2}^2
\nonumber
\\
&&
-a_6{\ovl a}_8z_{1}^2{\ovl z}_{1}{\ovl z}_{2}
-a_6{\ovl a}_{10}z_{1}^2{\ovl z}_{2}^2
+a_7{\ovl a}_8z_{1}z_{2}{\ovl z}_{1}^2
\nonumber
\\
&&
+a_7{\ovl a}_8z_{2}^2{\ovl z}_{1}{\ovl z}_{2}
-a_7{\ovl a}_9z_{1}z_{2}{\ovl z}_{2}^2
+a_8{\ovl a}_{10}z_{1}z_{2}{\ovl z}_{2}^2
)
\nonumber
\\
&&
+2(
a_3{\ovl a}_1z_{2}^2{\ovl z}_{1}^2
+a_4{\ovl a}_2z_{2}^2{\ovl z}_{1}^2
+a_7{\ovl a}_5z_{1}z_{2}{\ovl z}_{1}^2
-a_9{\ovl a}_5z_{2}^2{\ovl z}_{1}^2
\nonumber
\\
&&
-a_8{\ovl a}_6z_{1}z_{2}{\ovl z}_{1}^2
-a_{10}{\ovl a}_6z_{2}^2{\ovl z}_{1}^2
+a_8{\ovl a}_7z_{1}^2{\ovl z}_{1}{\ovl z}_{2}
\nonumber
\\
&&
+a_8{\ovl a}_7z_{1}z_{2}{\ovl z}_{2}^2
-a_9{\ovl a}_7z_{2}^2{\ovl z}_{1}{\ovl z}_{2}
+a_{10}{\ovl a}_8z_{2}^2{\ovl z}_{1}{\ovl z}_{2}
)
\nonumber
\\
&&
+\frac{4}{3}(
a_2{\ovl a}_3z_{1}^2{\ovl z}_{1}{\ovl z}_{2}
+a_2{\ovl a}_3z_{1}z_{2}{\ovl z}_{2}^2
+a_3{\ovl a}_2z_{1}z_{2}{\ovl z}_{1}^2
+a_3{\ovl a}_2z_{2}^2{\ovl z}_{1}{\ovl z}_{2}
)
\nonumber
\\
&&
+\frac{2}{3}(
a_2{\ovl a}_2z_{1}^2{\ovl z}_{1}^2
+a_3{\ovl a}_3z_{2}^2{\ovl z}_{2}^2
)
\nonumber
\\
&&
-\frac{15}{16}(
f_1^2z_{1}^2{\ovl z}_{1}^2
+f_4^2z_{2}^2{\ovl z}_{2}^2
)
-\frac{3}{4}(
f_1f_3z_{1}z_{2}{\ovl z}_{1}{\ovl z}_{2}
+f_2f_4z_{1}z_{2}{\ovl z}_{1}{\ovl z}_{2}
)
\nonumber
\\
&&
-\frac{5}{8}(
f_1f_2z_{1}^2{\ovl z}_{1}{\ovl z}_{2}
+f_1f_2z_{1}z_{2}{\ovl z}_{1}^2
+f_3f_4z_{1}z_{2}{\ovl z}_{2}^2
+f_3f_4z_{2}^2{\ovl z}_{1}{\ovl z}_{2}
)
\nonumber
\\
&&
-\frac{5f_2f_3}{24}(
z_{1}^2{\ovl z}_{1}{\ovl z}_{2}
+z_{1}z_{2}{\ovl z}_{1}^2
+z_{1}z_{2}{\ovl z}_{2}^2
+z_{2}^2{\ovl z}_{1}{\ovl z}_{2}
)
\nonumber
\\
&&
-\frac{1}{6}(
f_2^2z_{1}z_{2}{\ovl z}_{1}{\ovl z}_{2}
+f_3^2z_{1}z_{2}{\ovl z}_{1}{\ovl z}_{2}
)
\nonumber
\\
&&
-\frac{1}{8}(
f_2^2z_{1}^2{\ovl z}_{2}^2
+f_2^2z_{2}^2{\ovl z}_{1}^2
+f_3^2z_{1}^2{\ovl z}_{2}^2
+f_3^2z_{2}^2{\ovl z}_{1}^2
)
-\frac{5}{48}(
f_2^2z_{1}^2{\ovl z}_{1}^2
+f_3^2z_{2}^2{\ovl z}_{2}^2
)
\nonumber
\\
&&
+\frac{1}{16}(
f_1f_3z_{2}^2{\ovl z}_{1}^2
+f_1f_3z_{1}^2{\ovl z}_{2}^2
+f_2f_4z_{1}^2{\ovl z}_{2}^2
+f_2f_4z_{2}^2{\ovl z}_{1}^2
),
\label{cal-H4}
\end{eqnarray}
where $a_h$ ($h=1,\cdots , 10$) and $c_{\ell}$
($\ell = 1, \cdots , 13$) are the complex-valued
parameters chosen arbitrarily,
and $f_k$ ($k=1,\cdots,4$) the real-valued
parameters in ${\cal F}^{(3)}(q)$ (see (\ref{cal-F3})).
Namely, we have 46-degree-of freedom in the solution,
${\cal H}$, of the inverse problem of the PHOCP if
$(f_k)$ fixed.
Note that if $(a_h)$, $(c_{\ell})$ are chosen to be
\begin{equation}
\begin{array}{l}
\displaystyle{
a_1=a_3=a_5=a_8=a_9=0,
}
\\ \noalign{\vskip 6pt}
\displaystyle{
2a_2=2a_6=a_7=\frac{1}{4}, \; \; 
3a_4=a_{10}=\frac{3\mu}{8},
}
\end{array}
\label{special1-HH}
\end{equation}
and
\begin{equation}
\label{special2-HH}
c_{\ell}=0 \quad  (\ell=1,\cdots,13),
\end{equation}
respectively, and $(f_k)$ to be (\ref{f-HH}), ${\cal H}$
becomes equal to the one-parameter H{\' e}non-Heiles
Hamiltonian. After (\ref{cal-H})-(\ref{cal-H4}),
one might understand the necessity of computer algebra
in the inverse problem.
\subsection{The BDIC for PHOCPs and PHOQPs}
We wish to find the condition for $(a_h)$, $(c_{\ell})$
and $(f_k)$ to bring ${\cal H}$ into the PHOQP-Hamiltonian
\begin{equation}
\label{cal-F4}
{\cal F}^{(4)}(q,p)
=
\frac{1}{2}\sum_{j=1}^2 (p_j^2 + q_j^2)
+
(g_1q_1^4 + g_2q_1^3q_2 +g_3q_1^2q_2^2 
   +g_4 q_1q_2^3 +g_5q_2^4),
\end{equation}
where $(g_{\ell})$ are real-values parameters.
A straightforward calculation
shows the following (see \cite{Uwano2000} for detail) .
\begin{theorem}
\label{PHOQP-general}
A PHOCP-Hamiltonian ${\cal F}^{(3)}$
shares its BG-normal form with a certain PHOQP-Hamiltonian
${\cal F}^{(4)}$ up to degree-$4$ if and only if
the PHOCP-Hamiltonian ${\cal F}^{(3)}$ satisfies
the BDIC,
\begin{equation}
3(f_1f_3+f_2f_4)-(f_2^2+f_3^2)=0,
\label{BDIC-PHOCP}
\end{equation}
for PHOCPs.
Under (\ref{BDIC-PHOCP}), the PHOQP-Hamiltonian
${\cal F}^{(4)}$ sharing its BG-normal form with the
PHOCP-Hamiltonian ${\cal F}^{(3)}$ is equal to
\begin{eqnarray}
{\cal Q}&=&\frac{1}{2}\sum_{j=1}^2 (p_j^2 + q_j^2)
-\frac{5}{18}(9f_1^2+f_2^2)q_1^4
-\frac{10}{9}(3f_1+f_3)f_2 q_1^3q_2
-\frac{5}{3}(f_2^2+f_3^2)q_1^2q_2^2
\nonumber
\\
&\phantom{=}&
-\frac{10}{9}(3f_4+f_2)f_3 q_1q_2^3
-\frac{5}{18}(9f_4^2+f_3^2)q_2^4
,
\label{cal-Q}
\end{eqnarray}
where $(f_h)$ are subject to the BDIC (\ref{BDIC-PHOCP}).
\end{theorem}
We remark here that the BDIC (\ref{BDIC-PHOCP}) appears
again in Appendix.
\par
We are now in a position to show the integrability of the PHOQP
with ${\cal Q}$ with $(f_h)$ subject to (\ref{BDIC-PHOCP}).
It is easily seen that that ${\cal F}^{(4)}$ becomes ${\cal Q}$
under the substitution,
\begin{equation}
\begin{array}{l}
\displaystyle{
g_1= -\frac{5}{18}(9f_1^2+f_2^2),
\quad
g_2=-\frac{10}{9}(3f_1+f_3)f_2 ,
\quad
g_3=-\frac{5}{3}(f_2^2+f_3^2) ,
}
\\ \noalign{\vskip 6pt}
\displaystyle{
g_4=-\frac{10}{9}(3f_4+f_2)f_3 , 
\quad
g_5=-\frac{5}{18}(9f_4^2+f_3^2) .
}
\end{array}
\label{g-f}
\end{equation}
A long but straightforward calculation shows that
$(g_{\ell})$ given by (\ref{g-f}) with (\ref{BDIC-PHOCP})
satisfy the BDIC
\begin{equation}
\label{BDIC-PHOQP}
\begin{array}{l}
9g_2^2+4g_3^2-24g_1g_3-9g_2g_4=0,
\\ \noalign{\vskip4pt}
9g_4^2+4g_3^2-24g_3g_5-9g_2g_4=0,
\\ \noalign{\vskip4pt}
(g_2+g_4)g_3-6(g_1g_4+g_2g_5)=0 ,
\end{array}
\end{equation}
for PHOQPs (see Appendix). To summarize, we have the
following.
\begin{theorem}
\label{theorem-BDIC-CP-QP}
If a PHOCP and a PHOCP share the same BG-normal form
up to degree-$4$, then both oscillators are integrable
in the sense that they satisfy the BDIC.
\end{theorem}
\section{Concluding Remarks}
We have described the procedure in ANFER and the application of
ANFER to the BDIC for the PHOCPs. In the following,
Several remarks are made on ANFER itself and its application
\begin{description}
\item[{(1)}] (Improving the program):\quad
As pointed out at the end of section~4, the procedure
in ANFER based on the composition of canonical transformation
would save save the required-memory comparing with the procedure
given in section~2.3. There would be much room for improvement
on this direction. It would be done in writing-up process,
so that the author, not so familiar with writing programs,
is trying to improve ANFER with several collaborators
\cite{UwanoWEB}. Even if restricted to for PHOCPs,
the inverse problem is expected to have rich related subjects
listed below. Hence,those who are interested on improving
ANFER will be welcome writing-up, because the improvement on ANFER
is expected to advance various studies related with
the inverse problem of the BG-normalization.  
\item[{(2)}] (the separability): \quad
By a simple calculation, we see that the PHOCP-Hamiltonian
with
$
f_1=f_3/3=a+b, \; f_2/3=f_4=a-b
$
and the PHOQP-Hamiltonian with
$
g_1=g_3/6=g_5=-5(a^2+b^2), \; g_2=g_4=20(a^2+b^2)
$
share the same BG-normal form up to degree-$4$, where $a$ and $b$
are real-valued parameters \cite{Uwano2000}.
These oscillator Hamiltonians are well known to be separable 
in $q_1 \pm q_2$ (see \cite{Perelomov}),
Theorem~*** is hence understood to provide a significant
relation between the separability of the PHOCP-Hamiltonians
and that of PHOQP.
Since ${\cal F}^{(3)}$ with (\ref{BDIC-P3-a}) and ${\cal Q}$
are thought to include several classes of Hamiltonians separable in
several coordinate systems other than $q_1 \pm q_2$,
the separability will be worth studying extensively
from the BG-normalization viewpoint in future.
\item[{(3)}] (quantum bifurcation): \quad
Since the perturbed oscillators referred to in
theorems~\ref{PHOQP-general} and \ref{theorem-BDIC-CP-QP}
are integrable, their quantum spectra are expected to be
obtained exactly. These oscillators are hence expected to
provide good examples of the quantum bifurcation in the
BG-normalized Hamiltonian systems
\cite{Uwano1994,1995,1998,1999} to study
whether or not the quantum bifurcation
in the BG-normalized Hamiltonian system for those oscillators
approximates the bifurcation in these oscillators in a good extent.
\item[{(4)}] (integrability): \quad
As is easily seen,  the solution (\ref{cal-H}, \ref{cal-H3}, \ref{cal-H4})
of the inverse problem for ${\cal G}$ admits fifty real-valued
parameters. We may hence expect to obtain
other integrable systems, so-called the electromagnetic
type \cite{Hietarinta}, for example.
\end{description}
On closing this section, we wish to mention again
of the role of computer algebra in the inverse problem:
Without computer algebra, for example, in section~5,
it would have been very difficult to find
theorems~\ref{PHOQP-general} and \ref{theorem-BDIC-CP-QP}.
\par\smallskip\noindent
{\bf {\large Acknowledgment}}
\newline
The author wishes to thank Drs. S.I.Vinitsky, V.A.Rostovtsev
and A.A.Gutsev at the Joint Institute for Nuclear Research,
Dubna, Russia, for discussions.
This work is partly supported by the Grant-in-Aid Exploratory
Research no.~11875022 from the Ministry of Education, Science and
Culture, Japan. 
\section{Appendix: the BDIC}
In this Appendix, we review the BDIC briefly.
We will restrict our attention here to the two-degree-of-freedom
natural Hamiltonian systems associated with the Euclidean metric
because the PHOCPs and the PHOQPs dealt with in section~5 are of
such type.
\par
The theorem due to Bertrand and Darboux is stated as follows
(see \cite{Marshall,Yamaguchi}, for example).
\begin{theorem}[Bertrand-Darboux]
\label{BDIC-theorem}
Let $F$ be a natural Hamiltonian of the form,
\begin{equation}
\label{F}
F(q,p)
=
\frac{1}{2}\sum_{j=1}^{2}p_j^2 +V(q),
\end{equation}
where $V(q)$ a differentiable function in $q$.
Then, the following three statements are equivalent for
the Hamiltonian system with $F$. 
\begin{description}
\item[{(1)}]
There exists a set of real-valued constants,
$
(\alpha, \beta, \beta^{\prime}, \gamma, \gamma^{\prime})
\ne
(0,0,0,0,0)
$,
for which $V(q)$ satisfies
\begin{equation}
\label{BDIC-general}
\begin{array}{l}
\displaystyle{
\left(
\frac{\partial^2 V}{\partial q_2^2}
-\frac{\partial^2 V}{\partial q_1^2}
\right)
(-2\alpha q_1q_2-\beta^{\prime}q_2 -\beta q_1 + \gamma)
}
\\ \noalign{\vskip 6pt}
\displaystyle{
\phantom{xxxx}
+2 \frac{\partial^2 V}{\partial q_1 \partial q_2}
(\alpha q_2^2 -\gamma q_1^2 +\beta q_2 -\beta^{\prime}q_1
+\gamma^{\prime})
}
\\ \noalign{\vskip 6pt}
\displaystyle{
\phantom{xxxxxxxx}
+\frac{\partial V}{\partial q_1}
(6\alpha q_2 + 3\beta)
-
\frac{\partial V}{\partial q_2}
(6\alpha q_1 + 3\beta^{\prime})=0 .
}
\end{array}
\end{equation}
\item[{(2)}]
The Hamiltonian system with $F$ admits an integral of motion
quadratic in momenta.
\item[{(3)}]
The Hamiltonian $F$ is separable in either Cartesian, polar,
parabolic or  elliptic coordinates.
\end{description}
\end{theorem}
\indent
Due to the statement (2) in theorem~\ref{BDIC-theorem},
a natural Hamiltonian system with $F$ is always
integrable if (\ref{BDIC-general}) holds true.
In this regard, we refer to (\ref{BDIC-general})
as the Bertrand-Darboux integrability condition (BDIC).
\par
For the PHOCPs and the PHOQPs, Yamaguchi and Nambu
\cite{Yamaguchi} have given a more explicit expression
of the BDIC (\ref{BDIC-general}) convenient for section~5:
\begin{lemma}
\label{lemma-YN}
Let ${\cal F}^{(k)}(q,p)$ ($k=3,4$) be the Hamiltonians
of the form (\ref{cal-F3}) and (\ref{cal-F4}).
\begin{description}
\item[{(1)}]
For the PHOCP with ${\cal F}^{(3)}(q,p)$,
the BDIC (\ref{BDIC-general}) is equivalent to either of
the following conditions, (\ref{BDIC-P3-a}), 
(\ref{BDIC-P3-b}), or (\ref{BDIC-P3-c});
\begin{eqnarray}
&&
3(f_1f_3+f_2f_4)-(f_2^2+f_3^2)=0,
\label{BDIC-P3-a}
\\
&&
f_1=2f_3, \quad f_2=f_4=0,
\label{BDIC-P3-b}
\\
&&
f_4=2f_2, \quad f_1=f_3=0.
\label{BDIC-P3-c}
\end{eqnarray}
\item[{(2)}]
For the PHOQP with ${\cal F}^{(4)}(q,p)$,
the BDIC (\ref{BDIC-general}) is equivalent to either of
the following conditions, (\ref{BDIC-P4-a}) or
(\ref{BDIC-P4-b});
\begin{eqnarray}
&&
\left\{
\begin{array}{l}
g_3=2g_1=2g_5,
\\
g_2=g_4=0,
\end{array}
\right. 
\label{BDIC-P4-a}
\\
&&
\left\{
\begin{array}{l}
9g_2^2+4g_3^2-24g_1g_3-9g_2g_4=0,
\\ \noalign{\vskip4pt}
9g_4^2+4g_3^2-24g_3g_5-9g_2g_4=0,
\\ \noalign{\vskip4pt}
(g_2+g_4)g_3-6(g_1g_4+g_2g_5)=0.
\end{array}
\right.
\label{BDIC-P4-b}
\end{eqnarray}
\end{description}
\end{lemma}
The equations, (\ref{BDIC-P3-a}) and (\ref{BDIC-P4-b}),
are referred to as the BDIC (\ref{BDIC-PHOCP}) for PHOCPs and
(\ref{BDIC-PHOQP}) for PHOQPs, respectively, in section~5.
\par
In closing Appendix, we wish to mention little more
of the BDIC:
As stated in Theorem~\ref{BDIC-theorem},
the BDIC (\ref{BDIC-general}) provide not only a necessary
and sufficient condition for an existence of first integrals
quadratic in momenta \cite{Darboux,Whittaker} but also
for the separability in either Cartesian, polar, 
parabolic or elliptic coordinates  \cite{Marshall}. Indeed,
the BDIC has been studied repeatedly from various
viewpoints; the separation of variables
\cite{Marshall,Grosche}, the complete integrability
\cite{Perelomov}, so-called the direct method
\cite{Hietarinta} and the renormalization of Hamiltonian
equation \cite{Yamaguchi},
for example
(see also \cite{Whittaker} as an older reference and
the references in the above-cited literature).
\par


\begin{thebibliography}{99}
\bibitem{Uwano2000}
Y.~Uwano, J.~Phys.~{\bf A} {\bf 33}, 6635-6653 (2000).
%
\bibitem{CASC}
Y.~Uwano, N.~A.~Chekanov, V.~A.~Rostovtsev and S.~I.~Vinitsky,
{\it Computer Algebra in Scientific Computing}
ed Ganzha~V~G et al (Berlin: Springer-Verlag), p~441 (1999).
\bibitem{CPC}
N.~A.~Chekanov, V.~A.~Rostovtsev, Y.~Uwano and S.~I.~Vinitsky,
{\it Comput.~Phys.~Commun.} {\bf 126} 47 (2000).
\bibitem{JadPhys}
N.~A.~Chekanov, M.~Hongo, V.~A.~Rostovtsev, Y.~Uwano and
S.~I.~Vinitsky,
{\it Physics of Atomic Nuclei} {\bf 61} 2029 (1998). 
\bibitem{Gustavson}
F.~G.~Gustavson, {\it Astronomical Journal} {\bf 71}, 670 (1966).
\bibitem{Kummer}
M.~Kummer, {\it Commun. Math. Phys.} {\bf 48}, 53 (1976).
\bibitem{Cushman}
R.~Cushman, {\it Proc.~R.~Soc.} {\bf A382} 361 (1982).
\bibitem{Moser}
J.~K.~Moser, {\it Lectures on Hamiltonian Systems} Memoirs of
A.M.S. {\bf 81} (Providence: A.M.S.) p~10 (1968).
\bibitem{Goldstein}
H.~Goldstein, {\it Classical Mechanics} 2nd~ed.,
 p~241 (Reading: Addison-Wesley, 1950).  
\bibitem{Arnold}
V.~I.~Arnold,
{\it Mathematical Methods of Classical Mechanics}, p~215 
(New York: Springer-Verlag, 1980).
\bibitem{UwanoWEB}
Y.~Uwano, http://yang.amp.i.kyoto-u.ac.jp/\~{}uwano/ (2000).
\bibitem{Spivak}
M.~Spivak, {\it Calculus on Manifolds} (New York: Benjamin),
p~41 (1965).
\bibitem{Marshall}
I.~Marshall and S.~Wojciechowski, {\it J.~Math.~Phys.} {\bf 29},
1338 (1988).
\bibitem{Yamaguchi}
Y.~Y.~Yamaguchi and Y.~Nambu, {\it Prog. Theor. Phys.} {\bf 100},
199 (1998).
\bibitem{Darboux}
G.~Darboux, {\it Archives Neerlandaises} (ii) {\bf 6} 371 (1901).
\bibitem{Grosche}
C.~Grosche, G.~S.~Pogosyan and A.~N.~Sissakian,
{\it Fortschr. Phys.} {\bf 43}, 453 (1995).
\bibitem{Perelomov}
A.~M.~Perelomov, {\it Integrable Systems of Classical Mechanics
and lie Algebras} vol.1 (Basel: Birkh{\" a}user-Verlag),
p~81 (1990).
\bibitem{Hietarinta}
J.~Hietarinta, {\it Phys. Rep.} {\bf 147}, 87 (1987).
\bibitem{Whittaker}
E.~T.~Whittaker, {\it A Treatise on the Analytical Dynamics
of Particles and Rigid Bodies} 4th ed. (Cambridge: Cambridge U.P.),
p~332 (1937).
\bibitem{Uwano1994}
Y.~Uwano, {\it J. Phys. Soc. Jpn} {\bf 63 Suppl.~A}, 31 (1994).
\bibitem{1995}
Y.~Uwano, {\it J. Phys.} {\bf A28}, 2041 (1995).
\bibitem{1998}
Y.~Uwano, {\it Int. J. of Bifurcation and Chaos} {\bf 8}, 941 (1998).
\bibitem{1999}
Y.~Uwano, {\it Rep. Math. Phys.} {\bf 44}, 267 (1999).
\end{thebibliography}
\end{document}